\documentclass[11pt,a4paper]{article}

\usepackage{amsfonts}
\usepackage{amssymb}
\usepackage{amsthm}
\usepackage{epsfig}
\usepackage{color}
\usepackage{amsmath}
\usepackage{amsfonts}
\usepackage[english]{babel}
\usepackage{fontenc}
\usepackage{euscript}
\usepackage{vmargin}
\usepackage{amssymb}
\usepackage{latexsym}
\usepackage{enumerate}
\usepackage{graphicx}
\usepackage{hyperref}
 \usepackage{tikz}
\usepackage[titletoc,toc,title]{appendix}

\def\rain{\to +\infty}

\def\convn{\hspace{0.25em}{\atop{\scriptstyle
n\rain}}\hspace{-2.5em} \longrightarrow \hspace{0.2cm}}

\def\N{{\rm I\kern-.20em N}}
\def\R{{\rm I\kern-.20em R}}
\def\indi{{1\kern-.20em\rm I}}

 \newtheorem{lem}{Lemma}[section]

\newtheorem{rem}{Remark}

\newtheorem{prop}{Proposition}[section]

\newtheorem{definition}{Definition}[section]
\newtheorem{ex}{Example}[section]

\newcommand{\pg}{\hspace{0.6cm}}

\newcommand{\bdem} {\begin{proof}}
\newcommand {\edem}{\hfill $\square$ \end {proof}}

\begin{document}
\baselineskip 15pt \setcounter{page}{1}
\title{\bf \Large  Clustering of high values in random fields}
\author{{\small Helena Ferreira\footnote{ helenaf@ubi.pt},\ \  Lu\'{i}sa Pereira\footnote{ lpereira@ubi.pt}, \ \ Ana Paula Martins\footnote{ amartins@ubi.pt} }\\
\\
{\small\it   Department of Mathematics, University of Beira Interior, Portugal}\\
}
 \maketitle
 \baselineskip 15pt

\begin{quote}
{\bf Abstract:} The asymptotic results that underlie applications of extreme random fields often assume that the variables are located on a regular discrete grid, identified with $\mathbb{Z}^2$, and that they satisfy stationarity and isotropy conditions.
Here we extend the existing theory, concerning the asymptotic behavior of the maximum and the extremal index, to non-stationary and anisotropic random fields, defined over discrete subsets of $\mathbb{R}^2$. We show that, under a suitable coordinatewise long range dependence condition, the maximum may be regarded as the maximum of an approximately independent sequence of submaxima, although there may be high local dependence leading to clustering of high values. Under restrictions on the local path behavior of high values, criteria are given for the existence and value of the spatial extremal index which plays a key role in determining the cluster sizes and quantifying the strength of dependence between exceedances of high levels.
The general theory is applied to the class of max-stable random fields, for which the extremal index is obtained as a function of well-known tail dependence measures found in the literature, leading to a simple estimation method for this parameter. The results are illustrated with non-stationary Gaussian and 1-dependent random fields. For the latter, a simulation and estimation study is performed.

{\bf Key Words:}\ \  Random field, max-stable process, extremal dependence, spatial extremal index


\end{quote}

\section{Introduction}

\pg Extremes of variables like wind, temperature and precipitation can  affect anybody at any place. The potential consequences include increases in severe windstorms, flooding, wildfires, crop failure, population displacements and increased  mortality. Apart from their direct impacts, these events will also have indirect effects such as increased costs for strengthening infrastructure or higher insurance premiums.
When the interest lies in the study of variables measured at specifically-located monitors, such as the variables mentioned above, as well as air pollution, soil porosity or hydraulic conductivity, among others, spatial modeling is necessary, so random fields constitute an active area of current research.

The treatment of spatial and temporal dependence in random fields has been influenced by the multivariate Gaussian model, where the dependence is characterized by the covariance structures. However, this model excludes all the situations of marginal distributions with heavier tails than the Gaussian distribution, leaving aside a huge set of problems related to rare events. Extreme Value Theory plays an important role in these situations.

A considerable amount of work has been done in extending results of Extreme Value Theory to random fields which have ${\mathbb{Z}}^2$ as their parameter space. Although their lack of easy separation of past and future, a general version of the classical Extreme Types Theorem was given and the existence of the extremal index shown, by replacing a single global dependence restriction by several assumptions, each dealing with one coordinate direction, for which past-future separation is considered (\cite{Leadbetter2}, \cite{Pereira2}, among others).
Under local restrictions on the oscillations of the values of the random field, Ferreira and Pereira (\cite{Ferreira2}) and Pereira and Ferreira (\cite{Pereira2}) compute the extremal index from the joint distribution of a finite number of variables.

In a random field with high local dependence, an exceedance is likely to have neighboring exceedances, resulting in a clustering of exceedances, which leads to a compounding of events in the limiting point process of exceedances (\cite{Ferreira3}). 

The aforementioned results assumed that the variables are located on a regular grid, identified with $\mathbb{Z}^2$, and sometimes that they satisfy stationarity and isotropy conditions. 
This is a big restriction for the majority of the applications since usually spatial data are not regularly spaced, stationary and dependence is anisotropic, due to the presence of a main direction of dependence.

In this paper we extend the existing theory, concerning the asymptotic behavior of the maximum, to non-stationary and anisotropic random fields, $\textbf{Z}_S=\{Z(x):x\in S\}$, where $S=\bigcup_{n\geq1}A_n$ and $A=\left\{A_n\right\}_{n\geq 1}$ is an increasing sequence of sets of isolated points of $\mathbb{R}^2$, subject to conditions on long range and local dependencies. We will assume, without loss of generality, that the variables $Z(x), \ x\in S$, have common distribution $F$, being $\bar{F}$ the corresponding survival function. We will denote the maximum and the minimum of $Z(x)$ over $B\subset S$ by $\bigvee_{x\in B} Z(x)$ and $\bigwedge_{x \in B} Z(x)$, respectively. More precisely, in Section 2 we define an asymptotically independence condition under which we prove that $\bigvee_{x\in A_n} Z(x)$, $n\geq 1$, behaves asymptotically as independent maxima over a family of disjoint subsets of $A_n$. 

The way spatial extreme events interact is also of interest in spatial statistics. For example, an unusually stormy day at a particular location may be followed by another one at the same or a neighboring location. This type of dependence among spatial extremes can be summarized through the spatial extremal index of the sequence $\textbf{Z}_A=\{Z(x):x\in A_n\}_{n\geq 1}$.

\begin{definition}
The sequence $\textbf{Z}_A$ has spatial extremal index $\theta_A$ if, for each $\tau>0$ and any sequence of real numbers $\left\{u_n{(\tau)}\right\}_{n\geq 1}$ satisfying
\begin{equation}
E\left(\displaystyle\sum_{x\in A_n} \indi_{\{Z(x)>u_n\}}\right)\convn \tau,
\end{equation}
where $\indi_A$ denotes the indicator function of the event $A$, it holds that
$$
\displaystyle\lim_{ n\rightarrow +\infty} P\left(\displaystyle\bigvee_{x\in A_n} Z(x)\leq u_n(\tau)\right)=\exp(-\theta_A\tau).
$$
\end{definition}
\vspace{0.2cm}

The extremal index of $\textbf{Z}_A$ is the key parameter to relate the limiting distributions of $\bigvee_{x\in A_n}Z(x)$ and $\bigvee_{x \in A_n}\widehat{Z}(x)$, where $\widehat{{\bold {Z}}}_A=\{\widehat{Z}(x), x\in A_n\}_{n\geq 1}$ is a sequence of independent and identically random variables having the same distribution function $F$ as each variable of the sequence $\textbf{Z}_A$. In fact, if $f(n)$ is the number of locations on $A_n$ and there exists a sequence of real numbers $\{u_n(\tau)\}_{n\geq 1}$ satisfying (1.1), then
$$
P\left(\displaystyle\bigvee_{x\in A_n} \widehat{Z}(x)\leq u_n(\tau)\right)=F^{f(n)}(u_n(\tau)) \convn \exp\left(-\displaystyle\lim_{ n\rightarrow +\infty} f(n) \overline{F}(u_n(\tau))\right)=e^{-\tau}
$$
and
$$
P\left(\displaystyle\bigvee_{x\in A_n} Z(x)\leq u_n(\tau)\right) \convn \exp\left(-\displaystyle\lim_{ n\rightarrow +\infty} \theta f(n) \overline{F}(u_n(\tau))\right)=e^{-\theta_A \tau}.
$$
So,
\begin{enumerate}
\item $\displaystyle\lim_{ n\rightarrow +\infty} P\left(\displaystyle\bigvee_{x\in A_n} Z(x)\leq u_n(\tau)\right)=\displaystyle\lim_{ n\rightarrow +\infty} P\left(\displaystyle\bigvee_{x\in B_n} \widehat{Z}(x)\leq u_n(\tau)\right),
$
with $\sharp B_n=\theta_A f(n)\leq \sharp A_n$, that is, for the sequence of real levels $\{u_n(\tau)\}_{n\geq 1}$, $\displaystyle\bigvee_{x\in A_n}Z(x)$ behaves asymptotically as the maximum of less than $f(n)$ independent variables.\\
\item $\displaystyle\lim_{ n\rightarrow +\infty} P(\displaystyle\bigvee_{x\in A_n} Z(x)\leq u_n(\tau))=\displaystyle\lim_{ n\rightarrow +\infty} P(\displaystyle\bigvee_{x\in A_n} \widehat{Z}(x)\leq v_n(\theta \tau)$ where $v_n(\theta \tau)\cong F^{-1}\left(1-\frac{\theta \tau}{f(n)}\right)>u_n(\tau)\cong F^{-1}\left(1-\frac{\tau}{f(n)}\right)$, that is, $\displaystyle\bigvee_{x \in A_n}Z(X)$ behaves asymptotically as the maximum of the same number of independent variables but relatively to a level higher than $u_n(\tau)$.
\end{enumerate}
\vspace{0.2cm}

We may then deduce that the limit of the sequence $\{c_n\equiv P\left(\bigvee_{x\in A_n} Z(x)\leq u_n(\tau)\right)\}_{n\geq 1}$ is the same as the one we would obtain when considering $\{\widehat{c}_n\equiv P\left(\bigvee_{x\in A_n} \widehat{Z}(x)\leq u_n(\tau)\right)\}_{n\geq 1}$, if in $c_n$ we replace the levels $u_n(\tau), \ n\geq 1$, by  "appropriately close" levels $v_n(\tau'), \ n\geq 1$, with $\tau'<\tau$, or if we consider a sequence $\{B_n\}_{n\geq 1}$ with $B_n\subset A_n$ and $\sharp B_n\sim \theta f(n)$, instead of sequence $\{A_n\}_{n\geq 1}$. 
This suggests that over $A_n$ and considering the levels $u_n(\tau)$ we should not obtain isolated exceedances of $Z(x)$ contrarily to $\widehat{Z}(x)$, and therefore in this situation they occur in clusters. Later on we will prove a result that reinforces this intuition.

We finish Section 2 with an existence criteria for the extremal index of $\textbf{Z}_A$.

Section 3 contains the theory surrounding the maximum and the extremal index of $\textbf{Z}_A$ under restrictions on its exceedance local path behavior, which allow clustering of high values. Surprisingly we obtain a simple method for computing the extremal index of sequence $\textbf{Z}_A$ as the limit of a sequence of tail dependence coefficients.

Section 4 is devoted to the application of the results to max-stable random fields. There we introduce the notions of local and regional extremal indices and relate them with $\theta_A$. Based upon these relations a simple estimator for $\theta_A$ is given and its performance is analyzed with an anisotropic and non-stationary 1-dependent max-stable random field. Conclusions are drawn in Section 5 and the proofs are collected in the appendices.

\section{Asymptotic spatial independence}


\pg In this section, we show that, under a suitable long range dependence condition, the maximum of random fields defined over discrete subsets of $\mathbb{R}^2$, may be regarded as the maximum of an approximately independent sequence of submaxima, even though there may be high local dependence leading to clustering of high values.

The results are obtained through an extension of the methodology in Ferreira and Pereira (\cite{Pereira2}), for extremes on a regular grid, relying on the novelty of irregularly occurring extremes in space. 

The dependence structure used here is a coordinatewise long range dependence condition, which restricts dependence by limiting
$
\mid P\left(\bigvee_{x\in C \cup D} Z(x)\leq u_n\right) - P\left(\bigvee_{x\in C} Z(x)\leq u_n\right)P\left(\bigvee_{x\in D} Z(x)\leq u_n\right)\mid
$
with the two index sets $C,D\subset A_n$ being "separated" from each other by a certain distance $l_n$ along each direction.
 
Throughout we shall say that the pair $\left(I,J\right)$ is in ${\cal{S}}(\pi_i(A_n),l_n)$ if $I\subset S$ and $J\subset S$ are subsets of consecutive values of $\pi_i(A_n)$ separated by at least $l_n$ values of $\pi_i(A_n)$, where $\pi_i$, $i=1,2$, denote the cartesian projections.
The cardinality of the sets $A_n$ and $\pi_i(A_n)$, $i=1,2$, will be denoted by $\sharp A_n=f(n)$, $\sharp \pi_i(A_n)=f_i(n)$, $i=1,2$ and we will assume that $f(n)\rightarrow +\infty$ as $n\rightarrow +\infty$.

\vspace{0.3cm}
\begin{definition}
Let $\{u_n\}_{n\geq 1}$ be a sequence of real numbers. If there exist sequences of positive integers $l=\left\{l_n\right\}_{n\geq 1}$ and $k=\left\{k_n\right\}_{n\geq 1}$ such that
\begin{equation}
l_n\convn +\infty,\ \  k_n\convn +\infty, \ \ k_nl_n\frac{f_i(n)}{f(n)}\convn 0, \ for \ each \ i=1,2,
\end{equation}
and $k_n^2\alpha\left(l_n,u_n\right)\convn 0$, with
$$\alpha\left(l_n,u_n\right)=\displaystyle \sup \left| P\left(\displaystyle\bigvee_{x\in C \cup D} Z(x)\leq u_n\right) - P\left(\displaystyle\bigvee_{x\in C} Z(x)\leq u_n\right)P\left( \displaystyle\bigvee_{x\in D} Z(x)\leq u_n\right)\right|,$$
where the supremum is taken over sets $C$ and $D$ such that, for each $i=1,2$, $C, D \in {\cal{S}}(\pi_i(A_n),l_n)$, we say that the sequence $\textbf{Z}_A$ satisfies condition $D(u_n,k_n,l_n)$.
\end{definition}

Under $D(u_n,k_n,l_n)$-condition we have asymptotic independence of maxima over disjoint sets of locations, as shown in the following result.

\begin{lem}\label{lemma2.1}
Suppose that the sequence $\textbf{Z}_A$ satisfies condition $D(u_n,k_n,l_n)$ for a sequence of real numbers $\{u_n\}_{n\geq 1}$ such that
\begin{equation}
\{f(n)\overline{F}(u_n)\}_{n\geq 1} \ \ is  \ bounded.
\end{equation}
If $I^{(s)}$, $s,t=1,...,k_n$, are disjoint subsets of $\pi_1(A_n)$ and $J^{(t)}$, $s,t=1,...,k_n$, are disjoint subsets of $\pi_2(A_n)$, and
$$B_n^{(s,t)}=\pi_2^{-1}(J^{(t)})\cap \pi_1^{-1}(I^{(s)})\cap A_n,\ \ \ \ \  s,t=1,...,k_n,$$
then, as $n\rightarrow +\infty$, we have
$$\left|P\left(\displaystyle\bigvee_{x\in \displaystyle\cup_{s,t}B_n^{(s,t)}} Z(x)\leq u_n\right) - \displaystyle\prod_{s,t}P\left(\displaystyle\bigvee_{x\in B_n^{(s,t)}} Z(x)\leq u_n\right)\right|\longrightarrow 0.$$
\end{lem}

The next result proves that asymptotically the distribution of the maximum of $\textbf{Z}$ over $A_n$, $n\geq 1$, coincides with the distribution of the maximum of $\textbf{Z}$ over a union of conveniently chosen disjoint subsets of $A_n$, whenever condition $D(u_n,k_n,l_n)$ holds. 

The underlying idea to obtain the asymptotic distribution of the maximum of $\textbf{Z}$ over $A_n$, $n\geq 1\geq 1$, is to subdivide $A_n$ into $k_n^2$ disjoint subsets, $B_n^{(s,t)}$, $s,t=1,...,k_n$, using the following construction method of the family ${\cal{B}}_n=\left\{B_n^{(s,t)}:s,t=1,\ldots,k_n\right\}$:
\begin{itemize}
\item build $\pi_1^{-1}(I^{(s)})\cap A_n$, $s=1,...k_n$, with $I^{(s)}$, $s=1,...,k_n$, abutting subsets of consecutive values of $\pi_1 (A_n)$, maximaly chosen  for the condition $$\displaystyle\sum_{x\in \pi_1^{-1}(I^{(s)})\bigcap A_n}P(Z(x)>u_n)\leq \frac{1}{k_n}\displaystyle\sum_{x\in A_n}P(Z(x)>u_n)$$;\\
\item for each $s=1,...k_n$, build $\pi_2^{-1}(J^{(s,t)})\cap\pi_1^{-1}(I^{(s)})\cap A_n$, $t=1,...,k_n$, with $J^{(s,t)}$, $t=1,...,k_n$, contiguous subsets of $\pi_2 (\pi_1^{-1}(I^{(s)})\cap A_n)$ and maximally chosen such that $$\displaystyle\sum_{x\in \pi_2^{-1}(J^{(s,t)})\bigcap \pi_1^{-1}(I^{(s)})\bigcap A_n}P(Z(x)>u_n)\leq \frac{1}{k_n^{2}}\displaystyle\sum_{x\in A_n}P(Z(x)>u_n).$$
\end{itemize}

Figure 1 llustrates one possible set of disjoint blocks $B_n^{(s,t)}$, $s,t=1,...,k_n$, with $\bigcup_{s,t}B_n^{(s,t)}=A_n$ constructed through the previous method, for a particular $A_n$, $n\in\mathbb N$.

\begin{figure}[htb!]

\begin{center}

\resizebox{8cm}{!}{%

\tikzstyle{my help lines}=[lightgray,thin,dashed]

\begin{tikzpicture}

\draw [style=my help lines] (-0.5,-0.5) grid (11,11);

\draw [->] (-1,0) -- (11.5,0) node [below right]{$\Pi_1(x)$};

\draw [->] (0,-1) -- (0,11.5) node [left]{$\Pi_2(x)$};

\foreach \Point in { (3.4,4.4), (2.5,7.5), (2.8,6.8), (4.4,1.4), (5.6,2.6), (6.3,7.3), (7.2,9.2), (2.6,1.1), (2.9,3.9), (5.7,4.7), (7.6,8.6), (8.4,5.4), (8.5,1.5), (9.3,3.3), (9.8,8.8), (5.6,8.6), (2.6,4.6), (8.2,2.2), (4.1,7.1), (5.8,2.8), (2.9,7.9), (3.5,1.5), (5.9,8.9), (4.8,2.8), (6.1,8.1), (7.1,7.9), (7.6,6.6), (7.8,7.8), (8.4,9.4), (8.7,5.7), (9.6,3.6), (9.2,8.2), (9.8,1.8), (0.5,1.0), (0.7,0.3), (1.1,8.3), (1.4,0.9), (0.2,3.2), (1.3,2.6), (0.1,5.5), (0.7,1.1), (0.9,2.3), (0.8,4.4), (1.3,7.3), (1.1,8.9), (1.0,2.1), (0.8,0.3), (0.2,0.2), (0.1,1.5), (1.5,0.2), (1.1,6.2), (0.4,8.4), (0.6,1.0), (6.5,0.5), (6.1,1.0), (6.7,3.2), (7.1,2.3), (7.3,3.4), (7.5,3.8), (2.3,8.1), (3.1,8.7), (0.7,3.6), (0.7,6.4), (0.6,6.1), (0.1,7.2), (0.3,7.9), (1.5,2.0), (1.2,1.5), (2.0,5.5), (1.7,4.9), (2.0,3.5), (1.7,4.0), (2.1,3.1), (4.0,0.5), (3.0,0.8), (3.5,0.2), (3.0,2.0), (3.5,3.0), (4.0,3.7), (5.0,4.0), (3.0,5.5), (3.5,6.5), (4.0,8.0), (5.0,9.0), (5.0,7.5), (5.5,6.0), (6.5,6.7), (5.0,5.8), (5.0,0.2), (7.5,5.5), (8.0,6.0), (8.2,4.5), (7.5,0.5), (8.0,1.0), (7.2,0.3), (9.0,0.7), (9.5,0.5), (9.3,0.9), (9.0,2.5), (9.2,1.1), (9.0,4.2), (9.7,5.0), (9.7,6.3), (8.8,7.0), (9.5,7.5), (9.7,6.9), (9.8,8.0), (4.0,5.5)
}{

    \node at \Point {\textbullet};

};




\draw[ultra thick](1.0,0)-- (1,9.0);

\draw[ultra thick](2.6,0)-- (2.6,9.0);

\draw[ultra thick](4.5,0)-- (4.5,9.5);

\draw[ultra thick](7.0,0)-- (7.0,9.5);

\draw[ultra thick](8.6,0)-- (8.6,9.5);

\draw[ultra thick](10,0)-- (10,9.0);

\draw[ultra thick](0,0.5)-- (1,0.5);

\draw[ultra thick](0,1.4)-- (1,1.4);

\draw[ultra thick](0,2.7)-- (1,2.7);

\draw[ultra thick](0,5.4)-- (1,5.4);

\draw[ultra thick](0,6.8)-- (1,6.8);

\draw[ultra thick](0,8.7)-- (1,8.7);

\draw[ultra thick](0,8.7)-- (1,8.7);

\draw[ultra thick](1,1.3)-- (2.6,1.3);

\draw[ultra thick](2.6,3.0)-- (1,3.0);

\draw[ultra thick](1,9.0)-- (2.6,9.0);

\draw[ultra thick](1,8.0)-- (2.6,8.0);

\draw[ultra thick](1,6.0)-- (2.6,6.0);

\draw[ultra thick](1,4.5)-- (2.6,4.5);

\draw[ultra thick](2.6,1.0)-- (4.5,1.0);

\draw[ultra thick](2.6,2.5)-- (4.5,2.5);

\draw[ultra thick](2.6,4.0)-- (4.5,4.0);

\draw[ultra thick](2.6,6.0)-- (4.5,6.0);

\draw[ultra thick](2.6,7.5)-- (4.5,7.5);

\draw[ultra thick](2.6,9.0)-- (4.5,9.0);

\draw[ultra thick](4.5,1.1)-- (7.0,1.1);

\draw[ultra thick](4.5,3.0)-- (7.0,3.0);

\draw[ultra thick](4.5,5.0)-- (7.0,5.0);

\draw[ultra thick](4.5,7.0)-- (7.0,7.0);

\draw[ultra thick](4.5,8.5)-- (7.0,8.5);

\draw[ultra thick](4.5,9.5)-- (7.0,9.5);

\draw[ultra thick](7.0,1.3)-- (8.6,1.3);

\draw[ultra thick](7.0,2.8)-- (8.6,2.8);

\draw[ultra thick](7.0,4.7)-- (8.6,4.7);

\draw[ultra thick](7.0,6.1)-- (8.6,6.1);

\draw[ultra thick](7.0,8.0)-- (8.6,8.0);

\draw[ultra thick](7.0,9.5)-- (8.6,9.5);

\draw[ultra thick](8.6,1.0)-- (10.0,1.0);

\draw[ultra thick](8.6,1.0)-- (10.0,1.0);

\draw[ultra thick](8.6,3.0)-- (10.0,3.0);

\draw[ultra thick](8.6,4.6)-- (10.0,4.6);

\draw[ultra thick](8.6,6.5)-- (10.0,6.5);

\draw[ultra thick](8.6,7.7)-- (10.0,7.7);

\draw[ultra thick](8.6,9.0)-- (10.0,9.0);

\end{tikzpicture}

}

\caption{Example of a set of disjoint blocks $B_n^{(s,t)}$, $s,t=1,...,k_n$, with $\bigcup_{s,t}B_n^{(s,t)}=A_n$, for a particular $A_n$, $n\in\mathbb N$.
}

\end{center}

\end{figure}
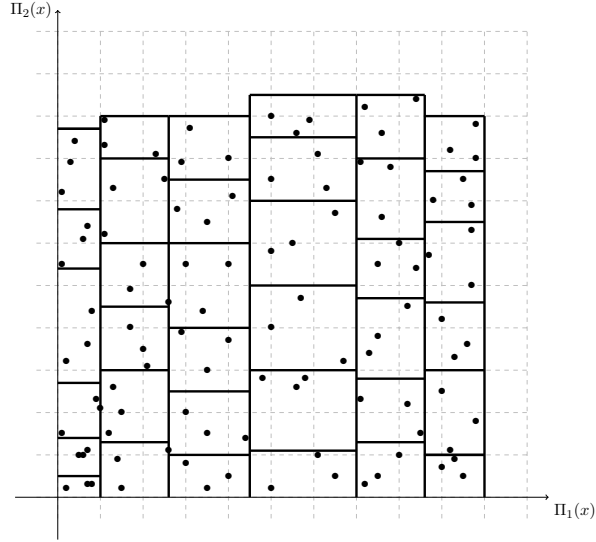

\vspace{0.3cm}

\begin{lem}\label{lemma2.2}
If the sequence ${\bold {Z}}_A$ satisfies condition $D(u_n,k_n,l_n)$, with $\{u_n\}_{n\geq 1}$  verifying (2.3) then, for each $n$, there exists a family  ${\cal{B}}_n$ of $k_n^2$ disjoint subsets of $A_n$, $B_n^{(s,t)}$, $s,t=1,...,k_n$, with $\sharp B_n^{(s,t)}\sim \frac{f(n)}{k_n^2}$ and such that
$$
P\left(\displaystyle\bigvee_{x\in A_n} Z(x)\leq u_n\right)-P\left(\displaystyle\bigvee_{x\in {\displaystyle\cup_{s,t}}B_n^{(s,t)}} Z(x)\leq u_n\right)\convn 0.
$$
\end{lem}

As a consequence of Lemmas \ref{lemma2.1} and \ref{lemma2.2} we can now state the following result concerning the asymptotic independence of maxima over $B_n^{(s,t)}, s,t=1,2,\ldots,k_n$.

\begin{prop}\label{proposition2.1}
If the sequence ${\bold {Z}}_A$ satisfies condition $D(u_n,k_n,l_n)$ with $\{u_n\}_{n\geq 1}$ verifying (2.3), then, for each $n$, there exists a family ${\cal{B}}_n$ of $k_n^2$ disjoint subsets of $A_n$, $B_n^{(s,t)}$, $s,t=1,\ldots,k_n$, with $\sharp B_n^{(s,t)}\sim \frac{f(n)}{k_n^2}$ and such that
$$
P\left(\displaystyle\bigvee_{x\in A_n} Z(x)\leq u_n\right)-\displaystyle\prod_{s,t}P\left(\displaystyle\bigvee_{x\in B_n^{(s,t)}} Z(x)\leq u_n\right) \convn 0.
$$
\end{prop}

\vspace{0.3cm}

\vspace{0.3cm}

The following result gives a convenient existence criteria for the extremal index of $\textbf{Z}_A$ and follows from Proposition 2.1: it depends on the local behavior of exceedances over $B_n^{(s,t)}$, $s,t=1,\ldots,k_n$, namely, on the limiting mean number of exceedances of $u_n$ by $\bigvee_{x\in B_n^{(s,t)}}Z(x)$, $s,t=1,\ldots,k_n$.

\begin{prop}\label{proposition2.2}
Suppose that the sequence ${\bold {Z}}_A$ satisfies condition $D(u_n(\tau),k_n,l_n)$, where $\{u_n(\tau)\}_{n\geq 1}$ is a sequence of real numbers satisfying (1.1) and ${\cal{B}}_n$ is a family of subsets of $A_n$ satisfying the conditions of Proposition \ref{proposition2.1}. Then, there exists the spatial extremal index, $\theta_A$, if and only if there exists
$$\displaystyle\lim_{ n\rightarrow +\infty} \displaystyle\sum_{B_n^{(s,t)}\in {\cal{B}}_n}P\left(\displaystyle\bigvee_{x\in B_n^{(s,t)}} Z(x)> u_n(\tau)\right),$$
and, in this case, we have
$$\theta_A=\displaystyle\lim_{ n\rightarrow +\infty} \frac{1}{f(n){\overline{F}}(u_n(\tau))}\displaystyle\sum_{B_n^{(s,t)}\in {\cal{B}}_n}P\left(\displaystyle\bigvee_{x\in B_n^{(s,t)}} Z(x)> u_n(\tau)\right).$$\\
\end{prop}

\vspace{0.3pt}

Next, we prove that the expected number of exceedances of the level $u_n(\tau)$ on the blocks $B_n^{(s,t)}$, $s,t=1,\ldots,k_n$ with at least one exceedance, converges to the reciprocal of the extremal index $\theta_A$.
We can verify that the greater the clustering tendency of high threshold exceedances (several exceedances on $B_n^{(s,t)}$) the smaller $\theta_A$ will be. For isolated exceedances of $u_n(\tau)$, we have $\theta_A=1$.

\begin{prop}\label{proposition2.3}
Suppose that the sequence ${\bold {Z}}_A$ satisfies condition $D(u_n(\tau),k_n,l_n)$, where $\{u_n(\tau)\}_{n\geq 1}$ is a sequence of real numbers satisfying (1.1) and ${\cal{B}}_n$ is a family of subsets of $A_n$ satisfying the conditions of Proposition \ref{proposition2.1}.
If the sequence $\textbf{Z}_A$ has spatial extremal index, $\theta_A$, then
$$\theta_A=\displaystyle\lim_{ n\rightarrow +\infty} \frac{1}{k_n^2}\displaystyle\sum_{B_n^{(s,t)}\in {\cal{B}}_n}
E^{-1}\left(\displaystyle\sum_{x\in B_n^{(s,t)}}\indi_{\{Z(x)>u_n(\tau)\}} \left|
\displaystyle\sum_{x\in B_n^{(s,t)}}\indi_{\{Z(x)>u_n(\tau)\}}>0\right)\right..$$

If
$$
\displaystyle\lim_{ n\rightarrow +\infty} E\left(\displaystyle\sum_{x\in B_n^{(s,t)}}\indi_{\{Z(x)>u_n(\tau)\}} \left|
\displaystyle\sum_{x\in B_n^{(s,t)}}\indi_{\{Z(x)>u_n(\tau)\}}>0\right)\right.=1,
$$
uniformly in $s,t\in \{1,\ldots,k_n\}$, then $\theta_A=1$.
\end{prop}

\section{Local spatial dependence}




\pg The asymptotic behavior of the maximum of non-stationary and anisotropic random fields, defined over discrete subsets of $\mathbb{R}^2$, subject to restrictions on the local path behavior of high values is now analyzed.

Criteria are given for the existence and value of the spatial extremal index, which plays a key role in determining the cluster sizes and quantifying the strength of dependence between exceedances of high levels.
To attain this goal, we first introduce a condition for modeling local mild oscillations of the random field. This condition is an extension to random fields of the $D''(u_n)$-condition found in Leadbetter and Nandagopalan (\cite{Leadbetter3}).

Throughout this section ${\cal{B}}_n$ will denote a family of subsets of $A_n$ in the conditions of Proposition 2.1.



\begin{definition}
If $V(x)$ is a finite set of neighbors of a point $x\in A$ and ${\cal{V}}=\{V(x):x\in A\}$, then the sequence $\textbf{Z}_A$ verifies condition $D''(u_n,{\cal{B}}_n,\cal{V})$ if, as $ n\rightarrow +\infty$, 
$$
k_n^2\displaystyle\sup_{B_n^{(s,t)}\in {\cal{B}}_n}\displaystyle\sum_{x\in B_n^{(s,t)}}P\left(Z(x)>u_n\geq \displaystyle\bigvee_{y\in {V}(x)} Z(y),\displaystyle\bigvee_{y\in B_n^{(s,t)}-V(x)} Z(y)>u_n\right) \longrightarrow 0.
$$
\end{definition}

\vspace{0.3cm}

Although the choice of the family $\cal{V}$ of neighborhoods can be conditioned by the nature of the practical problems under study, here we will illustrate the modeling with a natural choice based on the cardinal directions. Therefore, in what follows the initials N, E, S, W will represent, respectively, the cardinal directions North, East, South and West. The family of neighborhoods of $x\in A$ along directions E and N, will be denoted by ${\cal{V}}_{E,N}^{p,q,r}$, $p,q,r\in \mathbb{Z} \ \wedge \ q\leq 1$, and defined as 
\begin{eqnarray*}
{\cal{V}}_{E,N}^{p,q,r}=\{V(x):x\in A\},
\end{eqnarray*}
where
\begin{eqnarray*}
V(x)&=&\left\{y:\left(a_1(\pi_1(x))\leq \pi_1(y)\leq a_p(\pi_1(x)) \ \ \wedge \ \ a_{-q}(\pi_2(x))< \pi_2(y)\leq a_r(\pi_2(x))\right)\right.\\
&\ \ \ \ & \vee \ \ \left(\pi_1(x)=\pi_1(y) \ \ \wedge \ \ a_1(\pi_2(x))\leq \pi_2(y)\leq a_r(\pi_2(x))\right),x\in A\}, 
\end{eqnarray*}
and, for each $z\in \pi_i(A_n)$, 
$$
\ldots,a_{-2}(\pi_i(z)),a_{-1}(\pi_i(z)), a_{0}(\pi_i(z))=\pi_i(z), a_{1}(\pi_i(z)),a_{2}(\pi_i(z)),\ldots\\
$$
are the points before and after $\pi_i(z)$, in ascending order.
\vspace{0.3pt}

In Figure 2 we find an illustration of a neighborhood $V(x)\in {\cal{V}}_{E,N}^{p,q,r}$, $p,q,r\in \mathbb{Z} \ \wedge \ q\leq 1$.

\begin{figure}[htb!]

\begin{center}

\resizebox{7cm}{!}{%

\tikzstyle{my help lines}=[lightgray,thin,dashed]

\begin{tikzpicture}

 \draw [style=my help lines] (-0.5,-0.5) grid (13,13);

\draw [->] (-1,0) -- (13.5,0) node [below right]{$\Pi_1(x)$};

\draw [->] (0,-1) -- (0,13.5) node [left]{$\Pi_2(x)$};

\foreach \Point in {(2,2.5), (3,4), (3.5,6), (5.5,6), (4,5.5), (9,9), (10,9.5), (9.5,7), (4.5,7.5), (8,10.5), (1.5,1.5), (3,2), (4.1,2.6), (3.5, 3), (2,3.5), (4.2,3.2), (5.3,3), (2.5,2.3),(6.5,2.3),(7.5,3.5), (8.2,4.7),(2.5,4.5),(2.5,5.6), (3.2,7), (4.7,4.2), (5.2,4.2), (6,5), (6.4,6.5), (7.5,5.5), (8.5,5.5), (8,6.5),(9,6.3),(8.5,7.5),(7.5,7.8),(6.1,7.7),(5.3,8),(8,8.7),(10,7.8),(6.2,8.4),(7,9),(7.6,9.2),(7.1,10),(9.5,8,5), (10,8.7),(11,10),(10.8,9.2),(10.2,11.2),(9.8,10.8),(9.7,11.7),(9.5,11),(9,11),(4,6),(4,6.8)}{

    \node at \Point {\textbullet};

};

\node at (4,4.5) {$\times$};

\node [above right] at (4,4.5) {${\boldsymbol{x}}$};

\node [above right] at (3.3,7.3) {${\boldsymbol{V(x)}}$};

\draw[ultra thick](3.6,5)-- (4.5,5);

\draw[ultra thick](4.5,5)-- (4.5,4.7);

\draw[ultra thick](4.5,4.7)-- (6.6,4.7);

\draw[ultra thick](6.6,4.7)-- (6.6,7.3);

\draw[ultra thick](3.6,5)-- (3.6,7.3);

\draw[ultra thick](3.6,7.3)-- (6.6,7.3);

\end{tikzpicture}

}

\caption{Example of a neighborhood $V(x)\in {\cal{V}}_{E,N}^{p,q,r}$, with $p=10, q=-1$ and $r=11$}
\end{center}
\end{figure}
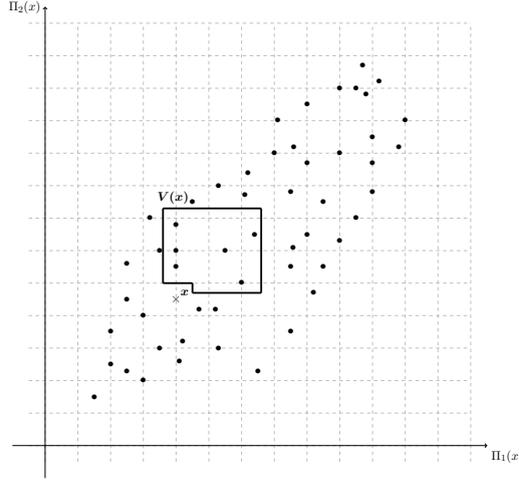

The following result proves that, asymptotically, the disjoint events that add to the probability of some exceedance of $u_n$ over $B_n^{(s,t)}$ are those where there occurs one exceedance of $u_n$ on the location $x\in B_n^{(s,t)}$ and the maximum on the $V(x)$ neighborhood is below $u_n$. That is, regarding these neighborhoods, $Z(x)$ is a local maximum.
\vspace{0.3cm}

\begin{lem}\label{lemma3.1}
Let $\{u_n\}_{n\geq 1}$ be a sequence of real numbers and suppose that sequence $\textbf{Z}_A$ satisfies condition $D''(u_n,{\cal{B}}_n,{\cal{V}})$. Then, for each $B_n^{(s,t)}\in {\cal{B}}_n$, we have
$$
P\left(\displaystyle\bigvee_{x\in B_n^{(s,t)}} Z(x)>u_n\right)=\displaystyle\sum_{x\in B_n^{(s,t)}}P\left(Z(x)>u_n>\displaystyle\bigvee_{y\in {V}(x)} Z(y)\right)+o\left(\frac{1}{k_n^2}\right).
$$
\end{lem}

\vspace{8pt}

As a consequence of Lemma \ref{lemma3.1} the extremal index of $\textbf{Z}_A$, $\theta_A$, can be viewed, as $n\rightarrow +\infty$, as the mean of tail dependence coefficients of the form
\begin{equation}
\lambda_n^{(\tau)}(V(x),x)=P\left(\left.\displaystyle\bigvee_{y\in V(x)} Z(y)>u_n(\tau)\right|Z(x)>u_n(\tau)\right),
\end{equation}
which are the tail dependence coefficient of Li (\cite{Li}). Observe also that if $\sharp V(x)=1$, then $\lambda_n^{(\tau)}(V(x),x)$ is the traditional upper tail dependence coefficient introduced far back in the sixties (Sibuya (\cite{Sibuya}), Tiago de Oliveira (\cite{Oliveira})).
\vspace{0.3cm}

\begin{prop}\label{proposition3.1}
Let $\{u_n(\tau)\}_{n\geq 1}$ be a sequence of real numbers satisfying (1.1).
If sequence $\textbf{Z}_A$ verifies conditions $D''(u_n(\tau),{\cal{B}}_n,{\cal{V}})$ and $D(u_n(\tau),k_n,l_n)$ then the spatial extremal index of $\textbf{Z}_A$, $\theta_A$, exists if and only if there exists
\begin{equation*}
\lambda_A=\displaystyle\lim_{n\rightarrow+\infty}\frac{1}{f(n)}\sum_{x\in A_n}\lambda_n^{(\tau)}(V(x),x),
\end{equation*}
and, in this case, we have
$$
\theta_A=1-\lambda_A.
$$
\end{prop}

\vspace{0.5cm}

Note that some models can verify condition $D''$ only for certain types of neighborhoods. Nevertheless, there exist models, as we shall see further on, that verify condition $D''(u_n,{\cal{B}}_n,{\cal{V}})$, for all ${\cal{V}}$. A particular case of such models are those that verify a local dependence restriction that leads to isolated exceedances, which we shall denominate condition $D'(u_n,{\cal{B}}_n)$ and define as follows:

\begin{definition}
The sequence $\textbf{Z}_A$ verifies condition $D'(u_n,{\cal{B}}_n)$ if, as $n\rightarrow \infty$,
\begin{equation*}
\displaystyle\sum_{B_n^{(s,t)}\in {\cal{B}}_n}\displaystyle\sum_{x,y\in B_n^{(s,t)}}P(Z(x)>u_n,Z(y)>u_n)\longrightarrow 0.
\end{equation*}
\end{definition}

\vspace{0.3cm}

This dependence condition, which bounds the probability of more than one exceedance of $u_n$ over a block $B_n^{(s,t)}$ with approximately $\frac{f(n)}{k_n^2}$ elements, together with condition $D(u_n,k_n,l_n)$  lead to an unit extremal index.
In fact, from Proposition \ref{proposition2.1} and condition $D'(u_n(\tau),{\cal{B}}_n)$ we have
\begin{eqnarray*}
\lefteqn{\displaystyle\lim_{ n\rightarrow +\infty}P\left(\displaystyle\bigvee_{x\in A_n} Z(x)\leq u_n(\tau)\right)}\\
&=&\exp\left(-\displaystyle\lim_{ n\rightarrow +\infty}\displaystyle\sum_{B_n^{(s,t)}\in {\cal{B}}_n}P\left(\displaystyle\bigvee_{x\in B_n^{(s,t)}} Z(x)> u_n(\tau)\right)\right)\\
&=&\exp\left(-\displaystyle\lim_{ n\rightarrow +\infty}\displaystyle\sum_{B_n^{(s,t)}\in {\cal{B}}_n}\displaystyle\sum_{x\in B_n^{(s,t)}}P\left(Z(x)> u_n(\tau))\right)+o(1)\right)\\
&=&\exp\left(-\displaystyle\lim_{ n\rightarrow +\infty}f(n)P\left(Z(x)> u_n(\tau)\right)\right)=\exp(-\tau),
\end{eqnarray*}
which proves the following result.

\begin{prop}\label{proposition3.2}
Let $\{u_n(\tau)\}_{n\geq 1}$ be a sequence of real numbers satisfying (1.1).
Suppose $\textbf{Z}_A$ satisfies conditions $D'(u_n(\tau),{\cal{B}}_n)$ and $D(u_n(\tau),k_n,l_n)$. Then, 
$$
P\left(\displaystyle\bigvee_{x\in A_n}Z(x)\leq u_n(\tau)\right)\convn \exp(-\tau).
$$
\end{prop}

\vspace{0.5cm}

We now present a class of Gaussian random fields that verifies the conditions established in Proposition \ref{proposition3.2}.
\vspace{0.5cm}

\begin{ex}
{\rm{
Let $\textbf{Z}_S=\{Z(x):x \in S\}$ be a standard Gaussian random field on $S\subset \mathbb{R}^2$ with correlations $r_{x,y}$, $x,y\in S$, such that
$$
\delta=\displaystyle\sup_{x,y\in S}\left|r_{x,y}\right|<1,
$$
where $S=\bigcup_{n\geq 1}A_n$ and $\{A_n\}_{n\geq 1}$ is an increasing sequence of sets of isolated points of $\mathbb{R}^2$, satisfying
$$
f(n)\convn +\infty, \ \ \ k_nl_n\frac{f_i(n)}{f(n)}\convn 0, \ \ \ \frac{f_i(n)}{f(n)}\convn 0, \ \ i=1,2,
$$
with $\{k_n\}_{n\geq 1}$ and $\{l_n\}_{n\geq 1}$ sequences of integer numbers verifying $k_n\rightarrow +\infty$ and $l_n\rightarrow +\infty$.
We will show that under the following correlation condition,
$$
\displaystyle\sup_{x,y \in A_n}\left|r_{x,y}\right|\leq r_n=o\left(\frac{1}{\log f(n)}\right),
$$
the sequence $\textbf{Z}_A=\{Z(x):x\in A_n\}_{n\geq 1}$ verifies conditions $D(u_n(\tau),k_n,l_n)$ and $D'(u_n(\tau),{\cal{B}}_n)$ with $u_n(\tau)=-\frac{\log \tau}{a_{f(n)}}+b_{f(n)}$, where $a_n=(2\log n)^{1/2}$ and $b_n=(2\log n)^{1/2}-\frac{1}{2}(2\log n)^{-1/2}(\log \log n+\log 4\pi)$.

By Corollary 4.2.9 of Leadbetter et al. (\cite{Leadbetter1}), we have
\begin{eqnarray*}
\lefteqn{\displaystyle \sum_{x,y \in B_n^{(s,t)}}P(Z(x)>u_n(\tau),Z(y)>u_n(\tau))}\\
&\leq&\displaystyle \sum_{x,y \in B_n^{(s,t)}}\left|P(Z(x)>u_n(\tau),Z(y)>u_n(\tau))-\Phi^2(u_n(\tau))\right|+\displaystyle \sum_{x,y \in B_n^{(s,t)}}\Phi^2(u_n(\tau))\\
&\leq&\displaystyle \sum_{x,y \in B_n^{(s,t)}}K\left|r_{x,y}\right|\exp\left(\frac{-u_n^2(\tau)}{1+\left|r_{x,y}\right|}\right)+\frac{f(n)}{k_n^2}\Phi^2(u_n(\tau)),\\
\end{eqnarray*}
where $\Phi$ denotes the standard Gaussian distribution function and $K$ is a constant depending on $\delta$.
\vspace{0.5cm}

Now, since $\displaystyle\sup_{x,y \in A_n}\left|r_{x,y}\right|\leq r_n=o\left(\frac{1}{\log f(n)}\right)$, by Lemma 4.3.2 of Leadbetter et al. (\cite{Leadbetter1}), we obtain
\begin{eqnarray*}
\lefteqn{\displaystyle \sum_{B_n^{(s,t)}}\displaystyle \sum_{x,y \in B_n^{(s,t)}}K\left|r_{x,y}\right|\exp\left(\frac{-u_n^2(\tau)}{1+\left|r_{x,y}\right|}\right)+k_n^2\frac{f(n)}{k_n^2}\Phi^2(u_n(\tau))}\\
&\leq&\displaystyle \sum_{x,y\in A_n}K\left|r_{x,y}\right|\exp\left(\frac{-u_n^2(\tau)}{1+\left|r_{x,y}\right|}\right)+o(1)\\
&=&o(1),\\
\end{eqnarray*}
proving that $\textbf{Z}_A$ verifies condition $D'(u_n(\tau),{\cal{B}}_n)$

Condition $D(u_n(\tau),k_n,l_n)$ follows from Corollary 4.2.4 of Leadbetter et al. (\cite{Leadbetter1}).
}}
\end{ex}

\vspace{0.3cm}

For other related results concerning Gaussian random fields we refer the readers to Piterbarg (\cite{Piterbarg}), Adler (\cite{Adler}), Berman (\cite{Berman}), Choi (\cite{Choi}) and Pereira (\cite{Pereira3}).

\vspace{0.3cm}

Proposition \ref{proposition3.1}, which states that the spatial extremal index $\theta_A$ is asymptotically equal to the mean of tail dependence coefficients is illustrated in the following example with a 1-dependent random field.

\begin{ex}
{\rm{
Let $\{Y(x):x\in \mathbb{Z}^2\}$ be an independent and identically distributed random field with common distribution function $F_Y(y)=\exp\left(\frac{-y^{-1}}{9}\right),\ y>0$, and define
\begin{equation}
Z(x)=\left\{
\begin{array}{lcl}
9Y(x) & , & \pi_1(x)=\pi_2(x)\\
\displaystyle\bigvee_{y\in E(x)}Y(y) & , &  \pi_1(x)\neq\pi_2(x),
\end{array}\right.
\end{equation}
where
$$
E(x)=\left\{y:a_{-1}(\pi_1(x))\leq \pi_1(y)\leq a_{1}(\pi_1(x))\wedge a_{-1}(\pi_2(x))\leq \pi_2(y)\leq a_{1}(\pi_2(x))\right\}.
$$

\begin{figure}[!htb]
\begin{center}
\includegraphics[scale=0.35]{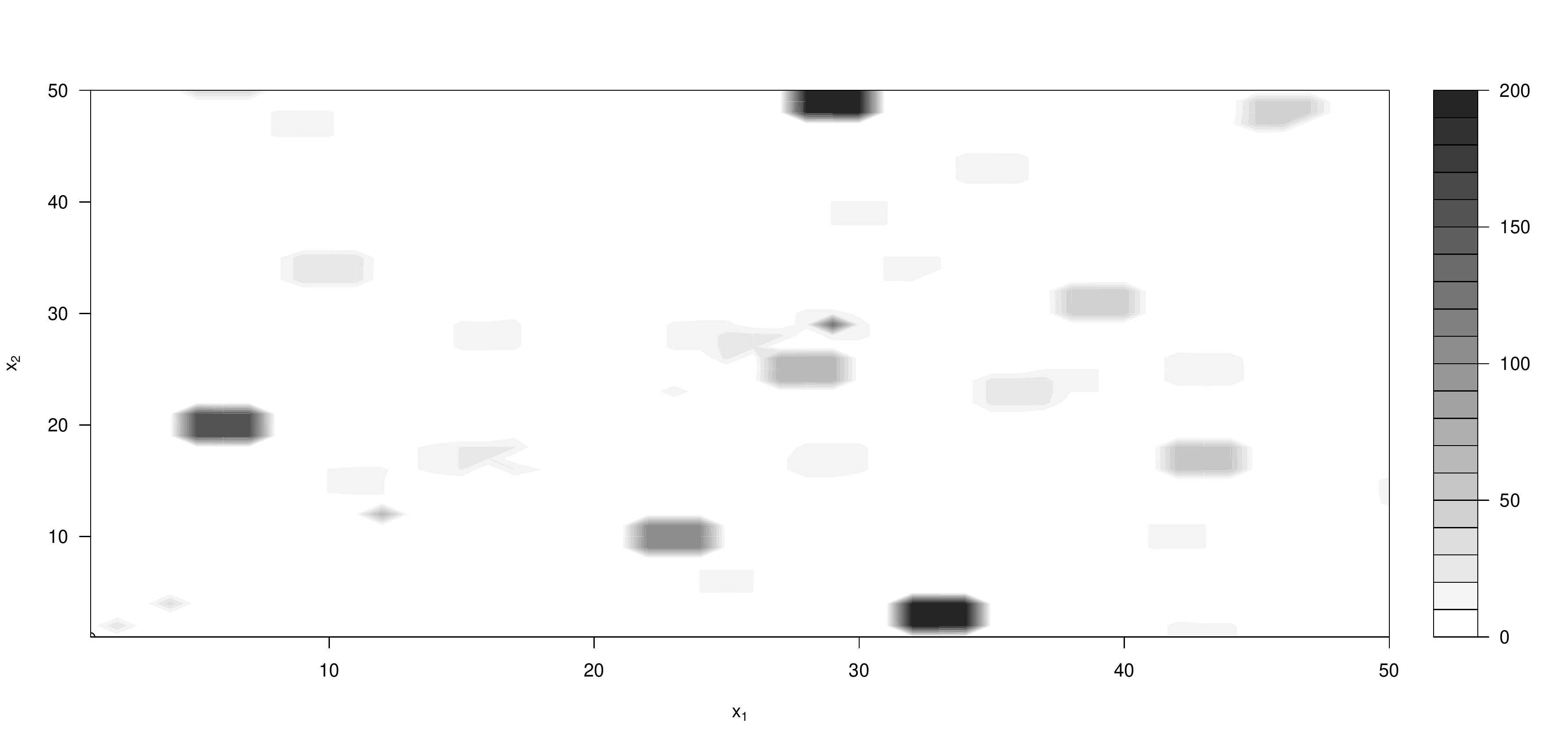}
\vspace{-0.5cm}\caption{\small Simulation of the random field $Z(x)$ defined in 3.5}
\label{fig:figura1}
\end{center}
\end{figure}
We will calculate the extremal index of sequence $\textbf{Z}_A=\{Z(x):x\in A_n\}_{n\geq 1}$, where
$$
A_n=\{x\in \mathbb{Z}^{2}:-n\leq \pi_1(x)\leq n,\ (-n)\vee (\pi_1(x)-n)\leq \pi_2(x)\leq (\pi_1(x)+n)\wedge n\}
$$
and
$$f(n)=3n^2+3n+1, \ \ f_i(n)=2n+1, \ i=1,2, \ \ \frac{f_i(n)}{f(n)}\convn 0, \ \ \bigcup_{n\geq 1}A_n=\mathbb{Z}^2,$$
by two different methods.

Note that the random field $\textbf{Z}=\{Z(x):x\in \bigcup_{n\geq 1}A_n\}$ is anisotropic and non-stationary with common unit Fréchet distribution, $F(y)=\exp(-y^{-1}),\ y>0$.

Let $\widehat{\textbf{Z}}=\{\widehat{Z}(x):x\in \mathbb{Z}^2\}$ be the associated random field of $\textbf{Z}$, id est, $\widehat{Z}(x)$, $x\in \mathbb{Z}^2$, are independent and identically distributed random variables having unit Fréchet distribution.

For a sequence of real numbers $\{u_n(\tau)\}_{n\geq 1}$ verifying (1.1), that is $u_n(\tau)=\frac{f(n)}{\tau}$, we have
\begin{equation*}
\lim_{n\rightarrow +\infty}P\left(\displaystyle\bigvee_{x\in A_n}\widehat{Z}(x)\leq u_n(\tau)\right)=\lim_{n\rightarrow +\infty}\left(\exp\left(-\frac{f(n)}{\tau}\right)^{-1}\right)^{f(n)}=\exp(-\tau).
\end{equation*}
On the other hand, for each $D_n=\{x\in A_n:\pi_1(x)=\pi_2(x)\}$ and $u_n(\tau)=\frac{f(n)}{\tau}$, it holds
\begin{equation*}
P\left(\displaystyle\bigvee_{x\in D_n}Z(x)>u_n(\tau)\right)\leq\sharp D_n \overline{F}(u_n(\tau))=\frac{2n+1}{f(n)}f(n)\overline{F}(u_n(\tau))\convn 0,
\end{equation*}
and consequently
\begin{eqnarray*}
P\left(\displaystyle\bigvee_{x\in A_n}Z(x)\leq u_n(\tau)\right)&=&P\left(\displaystyle\bigvee_{x\in A_n-D_n}Z(x)\leq u_n(\tau)\right)+o(1)\\
&=&P\left(\displaystyle\bigvee_{x\in A_{n+1}}Y(x)\leq u_n(\tau)\right)+o(1)\\
&=&(\exp(-u_n(\tau)^{-1}/9))^{3(n+1)^2+3(n+1)+1}\convn\exp\left(-\frac{\tau}{9}\right),\\
\end{eqnarray*}
so it holds $\theta_{A}=\frac{1}{9}$.

Let us now consider Proposition \ref{proposition3.1} for the computation of the extremal index of $\textbf{Z}_A$. 

Sequence $\textbf{Z}_A$ verifies condition $D(u_n(\tau),k_n,l_n)$, for all sequences of integer numbers $\{l_n\}_{n\geq 1}$ and $\{k_n\}_{n\geq 1}$ satisfying (2.2), since $\alpha_n(l_n,u_n(\tau))=0$ for $l_n\geq 2$.

Considering ${\cal{V}}={\cal{V}}_{S,W}^{1,1,0}$ the family of neighborhoods
\begin{eqnarray*}
V(x)=\{y:&& (a_{-1}(\pi_2(x))\leq \pi_2(y)<\pi_2(x)\wedge a_{-1}(\pi_1(x))\leq \pi_1(y)\leq\pi_1(x))\vee\\
&&(\pi_2(y)=\pi_2(x)\wedge a_{-1}(\pi_1(x))\leq \pi_1(y)<\pi_1(x))\},\\
\end{eqnarray*}
sequence $\textbf{Z}_A$ also verifies $D''(u_n(\tau),{\cal{B}}_n, {\cal{V}})$-condition since, for each $x\notin D_n$,
\begin{eqnarray*}
&&P\left(Z(x)>u_n(\tau)\geq\displaystyle\bigvee_{y\in V(x)}Z(y),\displaystyle\bigvee_{y\in B_n^{(s,t)}-V(x)}Z(y)>u_n(\tau)\right)\\
&\leq&P(Y(a)>u_n(\tau),Y(b)>u_n(\tau))=\overline{F}_Y^2(u_n(\tau)),\\
\end{eqnarray*}
for pairs of different locations $a$ and $b$, and consequently
\begin{eqnarray*}
&&k_n^2\displaystyle\sup_{B_n^{(s,t)}\in {\cal{B}}_n}\displaystyle\sum_{x\in B_n^{(s,t)}}P\left(Z(x)>u_n(\tau)\geq\displaystyle\bigvee_{y\in V(x)}Z(y),\displaystyle\bigvee_{y\in B_n^{(s,t)}-V(x)}Z(y)>u_n(\tau)\right)\\
&\leq&k_n^2\frac{f(n)}{k_n^2}\overline{F}_Y^2(u_n(\tau))+o(1)\\
&=&f(n)(1-\exp(-\tau/{9f(n)}))^2=o(1).\\
\end{eqnarray*}
We can then apply Proposition \ref{proposition3.1} to obtain
\begin{eqnarray*}
\theta_{A}&=&1-\lambda_A\\
&=&1-\displaystyle\lim_{n\to \infty}P\left(\displaystyle \bigvee_{y\in V(x)}Z(y)>u_n(\tau)\left| Z(x)>u_n(\tau)\right.\right)\\
&=&1-\displaystyle\lim_{n \to \infty}\displaystyle\sum_{x\in E(x)-\{(a_{-1}(\pi_1(x)),a_{-1}(\pi_2(x)))\}}\frac{P(Y(x)>u_n(\tau))}{P(Z(x)>u_n(\tau))}=\frac{1}{9}.
\end{eqnarray*}}}
\end{ex}

\section{Application to max-stable random fields}

\pg Max-stable random fields are very useful models for spatial extremes since under suitable conditions they are asymptotically  models for maxima of independent replications of random fields. Furthermore, all finite dimensional distributions of a max-stable process are multivariate extreme value distributions.

Within these random fields, it is important to identify dependence among extremes. In particular detecting asymptotic independence is fundamental and recently some authors have proposed measures of extreme dependence/independence with associated tests. With this in mind we compute, in this section, the extremal index of the class of max-stable random fields, as a function of well known extremal dependence coefficients found in literature, which will provide immediate estimators for $\theta_A$.

One convenient way to summarize the dependence structure of a max-stable random field $\textbf{Z}_S=\{Z(x):x\in S\}, \ S\subseteq \mathbb{R}^2$, with marginal distribution $F$, is through the extremal coefficient, $\displaystyle{\epsilon}_I$, of Schlather and Tawn (\cite{Schlather}), satisfying

$$P\left(\displaystyle\bigvee_{x\in I}Z(x)\leq y\right)=F^{{\displaystyle{\epsilon}}_I}(y), \ y\in\mathbb{R}, \ I\subseteq \mathbb{R}^2,
$$
which measures the extremal dependence between the variables indexed by the set $I$. Its simple interpretation as the effective number of independent variables indexed in $I$ from which the maximum is drawn has led to its use as a dependence measure in a range of practical applications. Another way to access the amount of extremal dependence of a random field is through a particular case of the tail dependence function introduced in Ferreira and Ferreira (\cite{Ferreira1}), defined as
\begin{equation}
\Lambda_U^{(I_1\mid I_2)}(y_1,y_2)=\lim_{t\to +\infty}P\left(\left.\displaystyle\bigvee_{x\in I_1}Z(x)>1-\frac{y_1}{t}\right| \displaystyle\bigvee_{x\in I_2}Z(x)>1-\frac{y_2}{t}\right), \ (y_1,y_2)\in \mathbb{R}_+^2,
\end{equation}
provided the limit exists. 

Note that for max-stable random fields, the limit given in (4.6) always exists.

The function $\Lambda_U^{(I_1\mid I_2)}(y_1,y_2)$ is a measure of the probability of occurring extreme values for the \linebreak maximum of the variables indexed in a region $I_1\subseteq \mathbb{R}_+^2$ given that the maximum of the variables indexed in another region $I_2$, with $I_1\cap I_2=\emptyset$, assumes an extreme value too.

At the unit point, we have
$$\Lambda_U^{(I_1\mid I_2)}(1,1)=\lim_{t\to +\infty}P\left(\left.\displaystyle\bigvee_{x\in I_1}Z(x)>1-\frac{1}{t}\right| \displaystyle\bigvee_{x\in I_2}Z(x)>1-\frac{1}{t}\right),
$$
which is related with the extremal coefficients of Schlatter and Tawn (\cite{Schlather}), $\displaystyle \epsilon$, in the following way
$$\Lambda_U^{(I_1\mid I_2)}(1,1)={\displaystyle{\epsilon}}_{I_1}+{\displaystyle{\epsilon}}_{I_2}-{\displaystyle{\epsilon}}_{I_1\cup I_2}.
$$\vspace{0.3cm}

The next result provides a connection between the dependence structure of the sequence of max-stable random fields $\textbf{Z}_A=\{Z(x):x\in A_n\}_{n\geq 1}$ and the limit of the sequence $\{\lambda_n^{(\tau)}(V(x),x)\}_{n\geq 1}$, with $\lambda_n^{(\tau)}(V(x),x)$ defined in (3.4).

We assume, without loss of generality that for each $x\in A_n$, $n\geq 1$, $Z(x)$ has a unit Fréchet distribution.

\begin{prop}\label{proposition4.1}
Let $\textbf{Z}_A=\{Z(x):x\in A_n\}_{n\geq 1}$ be a sequence of max-stable random fields with unit Fréchet margins. Then
$$\displaystyle\lim_{n\to +\infty}(1-\lambda_n^{(\tau)}(V(x),x))=\epsilon_{V(x)\bigcup\{x\}}-\epsilon_{V(x)}.
$$
\end{prop}

\vspace{0.3cm}

By combining the previous result with Proposition \ref{proposition3.1} we are able to conclude that if the extremal index of $\textbf{Z}_A$ exists, then it can be computed from $\epsilon_{V(x)\bigcup\{x\}}-\epsilon_{V(x)}$. We shall denote $\epsilon_{V(x)\bigcup\{x\}}-\epsilon_{V(x)}$ simply by $\theta_{\{x\}}$, $x\in A_n$, and name them local extremal indices.

If we consider the mean of local extremal indices for points on a region of $A_n$ we obtain a regional exremal index, formally defined as follows. 
\vspace{0.3cm}

\begin{definition}
Let $\textbf{Z}_A$ be a sequence of max-stable random fields with unit Fréchet margins. The extremal index of $\textbf{Z}_A$ over a region $R\subset A_n$, with $\sharp R<+\infty$, is defined as
$$
\theta_R=\frac{1}{\sharp R}\displaystyle\sum_{x\in R} \theta_{\{x\}}.
$$
\end{definition}

If $\textbf{Z}_A$ verifies conditions $D(u_n(\tau),k_n,l_n)$ and $D''(u_n(\tau),{\cal{B}}_n)$ then, for large $n$, its spatial extremal index, $\theta_A$, can be viewed as the mean of local extremal indices, as stated in the following result.

\begin{prop}\label{proposition4.2}
Suppose that the sequence of max-stable random fields $\textbf{Z}_A$ has marginal unit Fréchet distribution and verifies conditions $D(u_n(\tau),k_n,l_n)$ and $D''(u_n(\tau),{\cal{B}}_n)$. 
If
$$
1-\lambda_n^{(\tau)}(V(x),x)\convn \theta_{\{x\}}, \ uniformly \ in \ x,
$$
then
$$
\theta_A=\displaystyle\lim_{n\rightarrow +\infty}\frac{1}{f(n)}\displaystyle\sum_{x\in A_n}\theta_{\{x\}}.
$$
\end{prop}

\vspace{0.3cm}

We consider in the following example an anisotropic and non-stationary random field to illustrate Proposition \ref{proposition4.2}. 
\begin{ex}
{\rm{
Let $\{Y(x):x\in \mathbb{Z}^2\}$ be an independent and identically distributed random field with common distribution function $F_Y(y)=\exp\left(\frac{-y^{-1}}{3}\right),\ y>0$ and define the random field $\{Z(x), x\in (E\cup F)\times\mathbb{Z}\}$, with $E=\{4k:k\in\mathbb{Z}\}$ and $F=\{4k+3:k\in \mathbb{Z}\}$, in the following way
\begin{equation}
Z(x)=\left\{
\begin{array}{lcl}
\displaystyle\bigvee_{y\in U_1(x)} Y(y) & , & x\in E\times\mathbb{Z}\\
\displaystyle\bigvee_{y\in U_2(x)}Y(y) & , &  x\in F\times\mathbb{Z},
\end{array}\right.
\end{equation}
where
$$
U_1(x)=\{x,(a_1(\pi_1(x)),a_1(\pi_2(x))),(a_1(\pi_1(x)),a_{-1}(\pi_2(x)))\}
$$
and
$$
U_2(x)=\{(a_{-1}(\pi_1(x)),a_{-1}(\pi_2(x))),(a_1(\pi_1(x)),\pi_2(x)),(a_{-1}(\pi_1(x)),a_{1}(\pi_2(x)))\}.
$$

\begin{figure}[!htb]
\begin{center}
\includegraphics[scale=0.35]{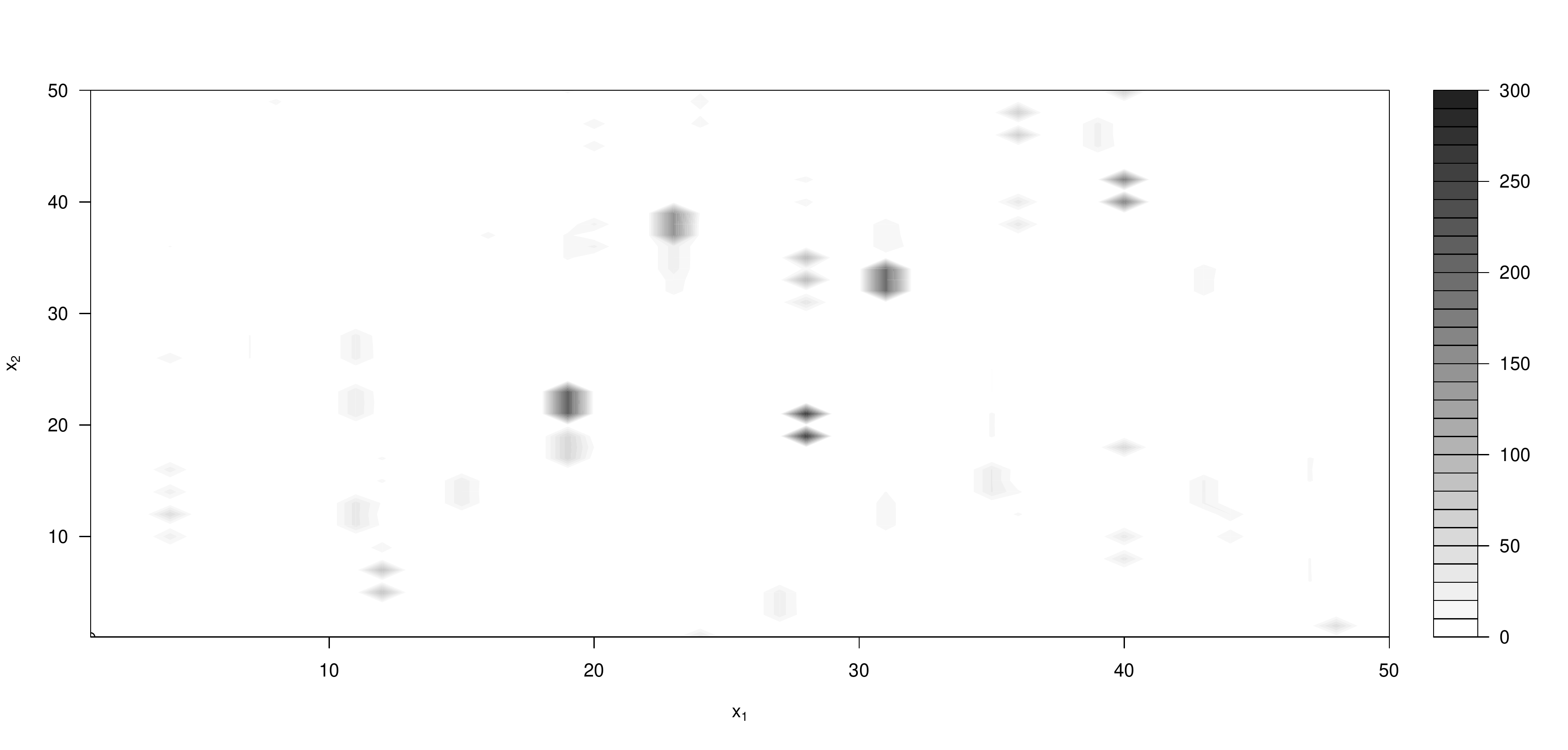}
\vspace{-0.5cm}\caption{\small Simulation of the random field $Z(x)$ defined in (4.7).}
\label{fig:figura1}
\end{center}
\end{figure}
The random field $\textbf{Z}=\{Z(x):x\in(E\cup F)\times \mathbb{Z}\}$ is anisotropic and non-stationary with common marginal distribution $F(x)=\exp(-x^{-1}), \ x\in(E\cup F)\times \mathbb{Z}$.

Let us consider 
$$
A_n=(E_n\cup F_n)\times G_n=(\{4k:-n\leq k\leq n\}\cup \{4k+3:-n\leq k\leq n\})\times\{-4n,\ldots,4n+3\},
$$
and consequently,
$$f(n)=8(2n+1)^2, \ \ f_1(n)=2(2n+1), \ \ f_2(n)=4(2n+1), \ \ \frac{f_i(n)}{f(n)}\convn 0, \bigcup_{n\geq 1}A_n=(E\cup F)\times \mathbb{Z}.$$

The sequence $\textbf{Z}_A=\{Z(x):x\in A_n\}_{n\geq 1}$ verifies condition $D(u_n,k_n,l_n)$, for all sequences of integer numbers $\{l_n\}_{n\geq 1}$ and $\{k_n\}_{n\geq 1}$ satisfying (2.2), since $\alpha_n(l_n,u_n)=0$ with $l_n\geq 2$, as well as condition $D''(u_n(\tau),{\cal{B}}_n, {\cal{V}})$ where $\{u_n(\tau)\}_{n\geq 1}$ is a sequence of real numbers satisfying (1.1) and $\cal{V}$ is the family of neighborhoods $V(x)=\{(\pi_1(x),a_1(\pi_2(x))), (\pi_1(x),a_2(\pi_2(x)))\}$. In fact, for each $x\in (E_n\times G_n)\cup (F_n\times G_n)$, we have
\begin{eqnarray*}
\lefteqn{P\left(Z(x)>u_n(\tau)\geq \displaystyle\bigvee_{y \in V(x)}Z(y),\displaystyle\bigvee_{y \in B_n^{(s,t)}-V(x)}Z(y)>u_n(\tau)\right)}\\
&=&P\left(Z(x)>u_n(\tau)\geq \displaystyle\bigvee_{y \in V(x)}Z(y),\displaystyle\bigvee_{y \in (B_n^{(s,t)}-V(x))\wedge\pi_2(y)>\pi_2(x)}Z(y)>u_n(\tau)\right)\\
\ \ &+&P\left(Z(x)>u_n(\tau)\geq \displaystyle\bigvee_{y \in V(x)}Z(y),\displaystyle\bigvee_{y \in (B_n^{(s,t)}-V(x))\wedge\pi_2(y)=\pi_2(x)}Z(y)>u_n(\tau)\right)\\
&\leq&\overline{F}_Z^2(u_n(\tau))\frac{f(n)}{k_n^2}+\overline{F}_Y^2(u_n(\tau))\frac{1}{k_n^2}.\\
\end{eqnarray*}
Therefore, for any family ${\cal{B}}_n$ in the conditions of the previously established results, we obtain
\begin{eqnarray*}
\lefteqn{k_n^2\sup_{B_n^{(s,t)}}\displaystyle\sum_{x\in B_n^{(s,t)}}P\left(Z(x)>u_n(\tau)\geq \displaystyle\bigvee_{y \in V(x)}Z(y),\displaystyle\bigvee_{y \in B_n^{(s,t)}-V(x)}Z(y)>u_n(\tau)\right)}\\
&\leq&\frac{(f(n)\overline{F}_Z(u_n(\tau)))^2}{k_n^2}+\frac{(f(n)\overline{F}_Y(u_n(\tau)))^2}{k_n^2}=o(1).\\
\end{eqnarray*}

From the definition of spatial extremal index of $\textbf{Z}_A$ we have $\theta_A=\frac{1}{2}$, since, for $u_n(\tau)=f(n)/\tau=8(2n+1)^2/\tau$, we obtain
\begin{eqnarray*}
P\left(\displaystyle\bigvee_{x \in A_n}Z(x)\leq u_n(\tau)\right)&=&P\left(\displaystyle\bigvee_{x \in E_n\times G_n}Z(x)\leq u_n(\tau)\right)P\left(\displaystyle\bigvee_{x \in F_n\times G_n}Z(x)\leq u_n(\tau)\right)\\
&=&(F_Y^3(u_n(\tau)))^\frac{f(n)}{2}\convn(\exp(-\tau))^\frac{1}{2}.\\
\end{eqnarray*}

On the other hand, from Proposition \ref{proposition4.2},
$$
\theta_A =1-\displaystyle\lim_{n\rightarrow\infty}\frac{1}{f(n)}\displaystyle\sum_{x\in A_n}\frac{P(Z(x)>u_n(\tau),Z(\pi_1(x),a_1(\pi_2(x)))>u_n(\tau) \vee Z(\pi_1(x),a_2(\pi_2(x)))>u_n(\tau))}{P(Z(x)>u_n(\tau))}.
$$
If $x\in E_n\times G_n$,
\begin{eqnarray*}
\lefteqn{\frac{P(Z(x)>u_n(\tau),Z(\pi_1(x),a_1(\pi_2(x)))>u_n(\tau) \vee Z(\pi_1(x),a_2(\pi_2(x)))>u_n(\tau))}{P(Z(x)>u_n(\tau))}}\\
&=&\frac{P(Z(x)>u_n(\tau),Z(\pi_1(x),a_1(\pi_2(x)))\leq u_n(\tau),Z(\pi_1(x),a_2(\pi_2(x)))>u_n(\tau))}{P(Z(x)>u_n(\tau))}\\
&=&\frac{P(Y(a_1(\pi_1(x)),a_1(\pi_2(x))>u_n(\tau))}{P(Z(x)>u_n(\tau))},\\
\end{eqnarray*}
otherwise, if $x\in F_n\times G_n$,
\begin{eqnarray*}
\lefteqn{\frac{P(Z(x)>u_n(\tau),Z(\pi_1(x),a_1(\pi_2(x)))>u_n(\tau) \vee Z(\pi_1(x),a_2(\pi_2(x)))>u_n(\tau))}{P(Z(x)>u_n(\tau))}}\\
&=&\frac{P(Z(x)>u_n(\tau),Z(\pi_1(x),a_1(\pi_2(x)))> u_n(\tau),Z(\pi_1(x),a_2(\pi_2(x)))>u_n(\tau))}{P(Z(x)>u_n(\tau))}+\\
&& \ \ \ \frac{P(Z(x)>u_n(\tau),Z(\pi_1(x),a_1(\pi_2(x)))> u_n(\tau),Z(\pi_1(x),a_2(\pi_2(x)))\leq u_n(\tau))}{P(Z(x)>u_n(\tau))}\\
&=&\frac{P(Y(a_{-1}(\pi_1(x)),a_1(\pi_2(x))>u_n(\tau))}{P(Z(x)>u_n(\tau))}+\frac{P(Y(a_{-1}(\pi_1(x)),\pi_2(x))>u_n(\tau))}{P(Z(x)>u_n(\tau))}.
\end{eqnarray*}
Thus,
\begin{eqnarray*}
\theta_A&=&1-\displaystyle\lim_{n\rightarrow\infty}\left(\frac{\overline{F}_Y(u_n(\tau))}{\overline{F}_Z(u_n(\tau))}\frac{(2n+1)\times 4(2n+1)}{f(n)}\right.+\\
&& \ \ \ \ \ \ \ \ \ \ \ \left.2\frac{\overline{F}_Y(u_n(\tau))}{\overline{F}_Z(u_n(\tau))}\frac{(2n+1)\times 4(2n+1)}{f(n)}\right)\\
&=&1-\displaystyle\lim_{n\rightarrow\infty}\left(\frac{f(n)\overline{F}_Y(u_n(\tau))}{f(n)\overline{F}_Z(u_n(\tau))}\times\frac{1}{2}+\frac{2f(n)\overline{F}_Y(u_n(\tau))}{f(n)\overline{F}_Z(u_n(\tau))}\times\frac{1}{2}\right)\\
&=&1-\left(\frac{\tau/3}{\tau}\times\frac{1}{2}+\frac{2\tau/3}{\tau}\times\frac{1}{2}\right)=\frac{1}{2},\\
\end{eqnarray*}}
as expected.}
\end{ex}

Although, in practical applications the conditions $D(u_n(\tau),k_n,l_n)$ and $D''(u_n(\tau),{\cal{B}}_n)$ are not easy to verify, the results of this section highlight the importance of $\theta_{\{x\}}$ in the study of locally occurring large observations in clusters. In a region $R\subseteq A_n$, the smaller the values of $\theta_{\{x\}}$, $x\in R$, the greater the propensity for clustering.

Note that, beyond the interpretation of the inverse proportionality between the value of $\theta_A$ and the propensity for clustering, small values of $\theta_{\{x\}}$ indicate a strong dependence between $Z(x)$ and $\{Z(y),y\in V(x)\}$.

\begin{subsection}{Estimation of the spatial extremal index}

\pg The spatial extremal index $\theta_A$ can be computed from the local extremal indices $\theta_{\{x\}}$, $x \in A_n$, and the latter are simply the diffrence between extremal coefficients at $V(x)\cup\{x\}$ and $V(x)$, as previously proved. Several estimators have already been studied in the literature (Krajina (\cite{Krajina}), Beirlant et al. (\cite{Beirlant}), Schlather and Tawn (\cite{Schlather}), among others) for extremal coefficients.
Ferreira and Ferreira (\cite{Ferreira1}) proposed a non-parametric estimator based on the following relation
$$
\epsilon_C=\frac{E(M(C))}{1-E(M(C))}, \ \ \ {\text {where}} \ \ \ M(C)=\bigvee_{x\in C}F(Z(x)).
$$
It considers sample means and is defined as
$$
\widehat{\epsilon}_{C}=\frac{\overline{M(C)}}{1-\overline{M(C})},
$$
where $\overline{M(C)}$ denotes the sample mean,
$$
\overline{M(C)}=\frac{1}{k}\sum_{i=1}^k \bigvee_{y \in C}\widehat{F}(Z^{(i)}(y))
$$
and $\widehat{F}$, is the (modified) empirical distribution function of $F$,
$$
\widehat{F}(u)=\frac{1}{k+1}\sum_{i=1}^k\indi_{\left\{Z^{(i)}(y)\leq u\right\}},
$$
where $Z^{(i)}(y)$, $i=1,\ldots,k$, are independent replications of $Z(y)$.

With this estimator of the extremal coefficient we propose the following estimator for the local extremal indices,
\begin{equation}
\widehat{\theta}_{\{x\}}=\widehat{\epsilon}_{\{x\}\cup V(x)}-\widehat{\epsilon}_{V(x)}  \ \ \ ,x\in A_n,
\end{equation}
which are consistent, given the consistency of the estimators $\widehat{\epsilon}_{\{x\}\cup V(x)}$ and $\widehat{\epsilon}_{V(x)}$, proved in Ferreira and Ferreira (\cite{Ferreira1}).

From Proposition \ref{proposition4.2} we know that, for large $n$, 
\begin{equation}
\theta_A\approx\frac{1}{f(n)}\displaystyle\sum_{x\in A_n}\theta_{\{x\}}.
\end{equation}
Therefore, if in (4.9) we replace $\theta_{\{x\}}$ with its estimator $\widehat{\theta}_{\{x\}}$ we obtain an estimator for the spatial extremal index $\theta_A$.

The finite sample behaviour of the estimator 
$$
\widehat{\theta}_A=\frac{1}{f(n)}\displaystyle\sum_{x\in A_n}\widehat{\theta}_{\{x\}}
$$
is analyzed on simulated data from the anisotropic and non-stationary random field $\{Z(x), x\in (E\cup F)\times\mathbb{Z}\}$ considered in Example 4.1. We simulated 10 times $k=100, 500$ and $1000$ independent random fields.

Table 1 shows the mean and mean square error (MSE) of the estimates for $k=100, 500, 1000$ and $f(n), \ n=1,10,20,30$.

\begin{table}

\begin{center}

\begin{tabular}{ccccccc} \hline

 & \multicolumn{6}{c}{k}\\ \cline{2-7}

& \multicolumn{2}{c}{100} & \multicolumn{2}{c}{500} & \multicolumn{2}{c}{1000}\\ \cline{2-7}

$f(n)$  &  {\scriptsize $E[\bullet]$ } & {\scriptsize $MSE[\bullet]$} &  {\scriptsize $E[\bullet]$ } & {\scriptsize $MSE[\bullet]$} & {\scriptsize $E[\bullet]$ } & {\scriptsize $MSE[\bullet]$} \\ \hline

f(1)=72  &  0.5253   & 7.6e-4 & 0.5140 & 2.2e-4 & 0.5099 & 1.1e-4 \\

f(10)=3528  &  0.5243  & 5.9e-4 & 0.5130 & 1.7e-4 & 0.5098 & 9.5e-5 \\

f(20)=13448& 0.5237 &5.6e-4
&0.5130 &1.7e-4&  0.5095& 9.0e-5\\

f(30)=29768& 0.5237& 5.6e-4& 0.5130& 1.7e-4&
0.5095& 9.1e-5\\\hline

\end{tabular}

\caption{Estimated mean values and mean square errors of estimator $\widehat{\theta}_A,$ for the random field of Example 4.1.}

\label{tab:tabela2}

\end{center}

\end{table}

As we can see from the values reported in Table 1, the estimator $\widehat{\theta}_A$ has quite a  good performance, with biases around 0.02 for small values of $k$ and around 0.01 for bigger values of $k.$  The values of $n$ considered have a small effect on the bias, nevertheless the variance decreases with $n$ and with $k.$ 


\end{subsection}

\section{Conclusion}

In this paper we establish existence criteria for the extremal index of a nonstationary and anisotropic random field, defined on $\mathbb{R}^2$.
Under restrictions on the local path behavior of exceedances, that allow clustering of high values, we obtain the extremal index as the limit of a sequence of upper tail dependence coefficients. For the particular case of max-stable random fields, we prove that the extremal index can be obtained as a function of extremal dependence coefficients. Based on this relation we give a simple estimator of the extremal index and we analyze its performance with an anisotropic and nonstationary 1-dependent random field. The simulation study results show the good performance of the proposed estimator.




\begin{appendix}\setcounter{equation}{0}

\section*{Appendix A: Proofs for Section 2}\setcounter{section}{0}

\refstepcounter{section}

\subsection*{Proof of Lemma \ref{lemma2.1}}

If all the subsets $I^{(s)}$ and $J^{(t)}$, $s,t=1,\ldots,k_n$, have less than $l_n$ elements, the result is trivial. On the other hand, if some $I^{(s_i)}$ has less than $l_n$ consecutive elements of $\pi_1(A_n)$, we can eliminate $B_n^{(s_i,t)}$ in the family $\left\{B_n^{(s,t)}:s,t=1,...k_n\right\}$, since
\begin{eqnarray*}
&&P\left(\displaystyle\bigvee_{x\in \bigcup_{s,t,s\neq s_i}B_n^{(s,t)}} Z(x)\leq u_n\right)-P\left(\displaystyle\bigvee_{x\in \bigcup_{s,t}B_n^{(s,t)}} Z(x)\leq u_n\right)\\
&=&P\left(\displaystyle\bigvee_{x\in\bigcup_t B_n^{(s_i,t)}}Z(x)>u_n\right)\\
&\leq&l_n \displaystyle\bigvee_{a} \frac{\sharp\left(\pi_1^{-1}(a)\cap A_n \right)}{f(n)} f(n) \overline{F}(u_n)=o\left(\frac{1}{k_n}\right)
\end{eqnarray*}
and
\begin{eqnarray*}
&&\left.\displaystyle \prod_{s,t,s\neq s_i} P\left(\displaystyle\bigvee_{x\in B_n^{(s,t)}} Z(x)\leq u_n\right)-\displaystyle\prod_{s,t}P\left(\displaystyle\bigvee_{x\in B_n^{(s,t)}} Z(x)\leq u_n\right)\right.\\
&\leq& 1-\displaystyle\prod_{t=1}^{k_n}P\left(\displaystyle\bigvee_{x\in B_n^{(s_i,t)}} Z(x)\leq u_n\right)\\
&\leq&\displaystyle\sum_{t=1}^{k_n}P\left(\displaystyle\bigvee_{x\in B_n^{(s_i,t)}} Z(x)> u_n\right)\\
&\leq& k_nl_n \displaystyle\bigvee_{a} \frac{\sharp\left(\pi_1^{-1}(a)\cap A_n \right)}{f(n)} f(n) \overline{F}(u_n)=o(1).
\end{eqnarray*}
With similar arguments, we can eliminate in the family $\left\{B_n^{(s,t)},s,t=1,...k_n\right\}$ the subsets $J^{(t_i)}$, $t_i\in\{1,\ldots,k_n\}$, that have less than $l_n$ elements.

To conclude, let us then assume that all the subsets $I^{(s)}$ and $J^{(t)}$, $s,t=1,2,\ldots,k_n$, have more than $l_n$ elements. Start by eliminating in each $I^{(s)}$ and $J^{(t)}$, respectively, the sets $I^{*{(s)}}$ and $J^{*^(t)}$ with the highest values $l_n$ values. The resulting sets $\overline{I}^{(s)}$, $s=1,\ldots,k_n$, belong to ${\cal{S}}(\pi_1(A_n),l_n)$, $s=1,\ldots,k_n$, and $\overline{J}^{(t)}$, $t=1,\ldots,k_n$,	belong to ${\cal{S}}(\pi_2(A_n),l_n)$. Now consider
$$
\overline{B}_n^{(s,t)}=\pi_1^{-1}(\overline{I}^{(s)})\cap \pi_2^{-1}(\overline{J}^{(t)})\cap A_n,
$$
and note that 
\begin{eqnarray*}
&&\left| P\left(\displaystyle\bigvee_{x\in {\displaystyle\cup_{s,t}}B_n^{(s,t)}} Z(x)\leq u_n\right) - \displaystyle\prod_{s,t} P\left(\displaystyle\bigvee_{x\in B_n^{(s,t)}} Z(x)\leq u_n\right)\right| \\
&\leq&\left| P\left(\displaystyle\bigvee_{x\in {\displaystyle\cup_{s,t}}B_n^{(s,t)}} Z(x)\leq u_n\right) - P\left(\displaystyle\bigvee_{x\in {\displaystyle\cup_{s,t}}\overline{B}^{(s,t)}} Z(x)\leq u_n\right)\right|\\
&+&\left| P\left(\displaystyle\bigvee_{x\in {\displaystyle\cup_{s,t}}\overline{B}^{(s,t)}} Z(x)\leq u_n\right) - \displaystyle\prod_{s,t} P\left(\displaystyle\bigvee_{x\in \overline{B}^{(s,t)}} Z(x)\leq u_n\right)\right| \\
&+&\left|\displaystyle\prod_{s,t} P\left(\displaystyle\bigvee_{x\in \overline{B}^{(s,t)}} Z(x)\leq u_n\right)-\displaystyle\prod_{s,t} P\left(\displaystyle\bigvee_{x\in B_n^{(s,t)}} Z(x)\leq u_n\right)\right|,\\
&=&2o(1)+k_n^2\alpha\left(l_n,u_n\right)=o(1),
\end{eqnarray*}
since
\begin{eqnarray*}
&& P\left(\displaystyle\bigvee_{x\in {\displaystyle\cup_{s,t}}\overline{B}_n^{(s,t)}} Z(x)\leq u_n\right) - P\left(\displaystyle\bigvee_{x\in {\displaystyle\cup_{s,t}}B_n^{(s,t)}} Z(x)\leq u_n\right)\\
&\leq&\displaystyle\bigvee_{i}2k_nl_n\frac{f_i(n)}{f(n)}f(n)\overline{F}(u_n)=o(1),
\end{eqnarray*}
condition $D(u_n,k_n,l_n)$ holds for $\textbf{Z}_A$ and
\begin{eqnarray*}
&&\displaystyle\prod_{s,t} P\left(\displaystyle\bigvee_{x\in \overline{B}_n^{(s,t)}} Z(x)\leq u_n\right)-\displaystyle\prod_{s,t}P\left(\displaystyle\bigvee_{x\in B_n^{(s,t)}} Z(x)\leq u_n\right)\\
&\leq&\displaystyle\sum_s\sum_t P\left(\displaystyle\bigvee_{x\in \overline{B}_n^{(s,t)}} Z(x)\leq u_n\right)-P\left(\displaystyle\bigvee_{x\in B_n^{(s,t)}} Z(x)\leq u_n\right)\\
&\leq&\displaystyle\sum_s\sum_t \left(l_n\sharp\pi_2(B_n^{(s,t)})+l_n\sharp\pi_1(B_n^{(s,t)})\right)\overline{F}(u_n)\\
&\leq&l_n \overline{F}(u_n)\left(k_n f_1(n)+k_n f_2(n)\right)\\
&\leq&\displaystyle\bigvee_{i}2k_nl_n\frac{f_i(n)}{f(n)}f(n)\overline{F}(u_n)=o(1).
\end{eqnarray*}


\subsection*{Proof of Lemma \ref{lemma2.2}}
First, we prove that the disjoint subsets $B_n^{(s,t)}$, $s,t=1,...,k_n$, constructed through the method mentioned before Lemma 2.2, has approximately $\frac{f(n)}{k_n^2}$ elements.
\begin{enumerate}
\item [1.]
In $\pi_1(A_n)$, let us consider the set $I^{(1)}$ of the first elements maximally constructed such that
$$
\displaystyle\sum_{x\in \pi_1^{-1}(I^{(1)})\bigcap A_n}P(Z(x)>u_n)\leq \frac{1}{k_n}\displaystyle\sum_{x\in A_n}P(Z(x)>u_n).
$$
In $\pi_1(A_n)-I^{(1)}$ let us consider the maximal set, $I^{(2)}$, of the first elements such that
$$
\displaystyle\sum_{x\in \pi_1^{-1}(I^{(2)})\bigcap A_n}P(Z(x)>u_n)\leq \frac{1}{k_n}\displaystyle\sum_{x\in A_n}P(Z(x)>u_n).
$$
Similarly, we obtain $k_n$ disjoint subsets of $A_n$,
$$
\pi_1^{-1}(I^{(1)})\cap A_n,\ldots,\pi_1^{-1}(I^{(k_n)})\cap A_n.
$$
\item [2.]
For each one of the previous subsets, let us consider an analogous decomposition using projection $\pi_2$.\\
In $\pi_2(A_n)$, we consider the set $J^{(s,1)}$ that consists of the first elements maximally chosen such that
$$
\displaystyle\sum_{x\in \pi_2^{-1}(J^{(s,1)})\bigcap \pi_1^{-1}(I^{(s)})\bigcap A_n}P(Z(x)>u_n)\leq \frac{1}{k_n^{2}}\displaystyle\sum_{x\in A_n}P(Z(x)>u_n).
$$
Using the same technique, we obtain the following $k_n$ subsets of $\pi_1^{-1}(I^{(s)})\cap A_n$,
$$
B_n^{(s,t)}=\pi_2^{-1}(J^{(s,t)})\cap \pi_1^{-1}(I^{(s)})\cap A_n, \ \ t=1,\ldots,k_n.
$$
\end{enumerate}
\vspace{0.3cm}
Next, we will prove that $\sharp B_n^{(s,t)}\sim\frac{f(n)}{k_n^2}.$
Since
$$
\displaystyle\sum_{x\in \pi_2^{-1}(J^{(s,t)})\bigcap \pi_1^{-1}(I^{(s)})\bigcap A_n}P(Z(x)>u_n)\leq \frac{1}{k_n^{2}}\displaystyle\sum_{x\in A_n}P(Z(x)>u_n)
$$
we obtain
$$
\sharp B_n^{(s,t)}\leq \frac{1}{k_n^2} f(n),
$$
and, from the maximality criteria used in the construction of $J^{(s,t)}$, it follows that
$$
\sharp B_n^{(s,t)}\overline{F}(u_n)+o\left(\frac{1}{k_nl_n}\right)>\frac{1}{k_n^2}f(n) \overline{F}(u_n).
$$
Therefore,
$$
1\geq \frac{\sharp B_n^{(s,t)}}{\frac{f(n)}{k_n^2}}\geq 1-\frac{k_n^2}{F(n)}o\left(\frac{1}{k_nl_n}\right),
$$
which allows us to conclude that $\sharp B_n^{(s,t)}\sim \frac{f(n)}{k_n^2}$, since $\frac{k_n^2}{f(n)}\convn 0$.\\
Now,
\begin{eqnarray*}
&&P\left(\displaystyle\bigvee_{x\in A_n} Z(x)\leq u_n\right)-P\left(\displaystyle\bigvee_{x\in {\displaystyle\cup_{s,t}}B_n^{(s,t)}} Z(x)\leq u_n\right)\\
&\leq&\displaystyle \sum_{x\in A_n-\displaystyle\cup_{s,t}B_n^{(s,t)}}P(Z(x)>u_n)\\
&=&\frac{\sharp\left( A_n-\displaystyle\bigcup_{s,t}B_n^{(s,t)}\right)}{f(n)}f(n)\overline{F}(u_n)
\end{eqnarray*}
and
\begin{eqnarray*}
&&\frac{f(n)-\sharp\displaystyle\bigcup_{s,t}B_n^{(s,t)}}{f(n)}=\frac{f(n)-\sharp\displaystyle\sum_{s,t}B_n^{(s,t)}}{f(n)}\\
&=&1-\frac{1}{k_n^2}\displaystyle\sum_{s,t}\frac{\sharp B_n^{(s,t)}}{\frac{f(n)}{k_n^2}}\leq 1-\displaystyle\bigwedge_{s,t}\frac{ \sharp B_n^{(s,t)}}{\frac{f(n)}{k_n^2}}=o(1),
\end{eqnarray*}
which concludes the proof.


\subsection*{Proof of Proposition \ref{proposition2.2}}

From Proposition \ref{proposition2.1}, we have
\begin{eqnarray*}
&&\displaystyle\lim_{ n\rightarrow +\infty}P\left(\displaystyle\bigvee_{x\in A_n} Z(x)\leq u_n(\tau)\right)\\
&=&\displaystyle\lim_{ n\rightarrow +\infty}\displaystyle\prod_{B_n^{(s,t)}\in {\cal{B}}_n}P\left(\displaystyle\bigvee_{x\in B_n^{(s,t)}} Z(x)\leq u_n(\tau)\right)\\
&=&\exp\left\{-\displaystyle\lim_{ n\rightarrow +\infty}\displaystyle\sum_{B_n^{(s,t)}\in {\cal{B}}_n}-\ln\left(1-P\left(\displaystyle\bigvee_{x\in B_n^{(s,t)}} Z(x)> u_n(\tau)\right)\right)\right\}\\
&=&\exp\left\{-\displaystyle\lim_{ n\rightarrow +\infty}\displaystyle\sum_{B_n^{(s,t)}\in {\cal{B}}_n} P\left(\displaystyle\bigvee_{x\in B_n^{(s,t)}} Z(x)> u_n(\tau)\right)\right\}.
\end{eqnarray*}
Then, if there exists $\lim_{ n\rightarrow +\infty} \sum_{B_n^{(s,t)}\in {\cal{B}}_n}P\left(\bigvee_{x\in B_n^{(s,t)}} Z(x)> u_n(\tau)\right)$, the result follows from the definition of $\theta_A$.


\subsection*{Proof of Proposition \ref{proposition2.3}}
Follows from Proposition \ref{proposition2.2}, since
\begin{eqnarray*}
\theta_A&=&\displaystyle\lim_{n\rightarrow +\infty} \frac{1}{f(n){\overline{F}}(u_n(\tau))}\displaystyle\sum_{B_n^{(s,t)}\in {\cal{B}}_n}P\left(\displaystyle\bigvee_{x\in B_n^{(s,t)}} Z(x)> u_n(\tau)\right)\\
&=&\displaystyle\lim_{ n\rightarrow +\infty} \frac{1}{k_n^2} \displaystyle\sum_{B_n^{(s,t)}\in {\cal{B}}_n}\frac{P\left(\displaystyle\bigvee_{x\in B_n^{(s,t)}} Z(x)> u_n(\tau)\right)}{\frac{f(n)}{k_n^2}{\overline{F}}(u_n(\tau))}\\
&=&\displaystyle\lim_{n\rightarrow +\infty} \frac{1}{k_n^2} \displaystyle\sum_{B_n^{(s,t)}\in {\cal{B}}_n}\frac{P\left(\displaystyle\sum_{x\in B_n^{(s,t)}}\indi_{\{Z(x)>u_n(\tau)\}}>0\right)}{E\left(\displaystyle\sum_{x\in B_n^{(s,t)}}\indi_{\{Z(x)>u_n(\tau)\}}\right)}.
\end{eqnarray*}


\end{appendix}

\begin{appendix}\setcounter{equation}{0}

\section*{Appendix B: Proofs for Section 3}\setcounter{section}{0}

\refstepcounter{section}

\subsection*{Proof of Lemma \ref{lemma3.1}}

We can write
\begin{eqnarray*}
&&P\left(\displaystyle\bigvee_{x\in B_n^{(s,t)}} Z(x)>u_n\right)\\
&=&P\left(\displaystyle\bigcup_{x\in B_n^{(s,t)}}\left\{Z(x)>u_n\geq \displaystyle\bigvee_{\{y\in B_n^{(s,t)}:\pi_1(x)>\pi_1(y)\vee(\pi_1(y)=\pi_1(x)\wedge \pi_2(y)>\pi_2(x))\}} Z(y)\right\}\right)\\
&=&\displaystyle\sum_{x\in B_n^{(s,t)}}P\left(Z(x)>u_n\geq \displaystyle\bigvee_{\{y\in B_n^{(s,t)}:\pi_1(x)>\pi_1(y)\vee(\pi_1(y)=\pi_1(x)\wedge \pi_2(y)>\pi_2(x))\}} Z(y)\right)\\
&=&\displaystyle\sum_{x\in B_n^{(s,t)}}P\left(Z(x)>u_n\geq \displaystyle\bigvee_{y\in V(x)} Z(y)\right)+o\left(\frac{1}{k_n^2}\right),
\end{eqnarray*}
since, by $D''(u_n,{\cal{B}}_n,{\cal{V}}_{E,N}^{p,q,r})$-condition, we have
\begin{eqnarray*}
&\left.\right.&\displaystyle\sum_{x\in B_n^{(s,t)}}P\left(Z(x)>u_n\geq \displaystyle\bigvee_{y\in V(x)} Z(y)\right)\\
&\left.\right.&\left.\left.\left.\right.\right.\right. -\displaystyle\sum_{x\in B_n^{(s,t)}}P\left(Z(x)>u_n\geq \displaystyle\bigvee_{\{y\in B_n^{(s,t)}:\pi_1(x)>\pi_1(y)\vee(\pi_1(y)=\pi_1(x)\wedge \pi_2(y)>\pi_2(x))\}} Z(y)\right)\\
&\left.\right.&\leq\displaystyle\sup_{B_n^{(s,t)}\in {\cal{B}}_n}\displaystyle\sum_{x\in B_n^{(s,t)}}P\left(Z(x)>u_n\geq \displaystyle\bigvee_{y\in V(x)} Z(y),\displaystyle\bigvee_{y\in B_n^{(s,t)}-V(x)} Z(y)>u_n\right)=o\left(\frac{1}{k_n^2}\right).
\end{eqnarray*}


\begin{rem}
Instead of the decomposition of $\{\bigvee_{x\in B_n^{(s,t)}} Z(x)>u_n\}$ as the union of events $\{Z(x)>u_n\}$, $x\in B_n^{(s,t)}$, where $x$ is the location with the highest coordinates, previously considered, we can consider other decompositions. For example, the events $\{Z(x)>u_n\}$ where $x\in B_n^{(s,t)}$ is the location with the lowest abscissa and biggest ordinate, lead to the following decomposition
\begin{eqnarray*}
\left\{\displaystyle\bigvee_{x\in B_n^{(s,t)}} Z(x)>u_n\right\}=\displaystyle\bigcup_{x\in B_n^{(s,t)}}\left\{Z(x)>u_n\geq \displaystyle\bigvee_{\{y\in B_n^{(s,t)}:\pi_1(y)<\pi_1(x)\vee(\pi_1(y)=\pi_1(x)\wedge \pi_2(y)>\pi_2(x))\}} Z(y)\right\}.
\end{eqnarray*}
Now, denoting by ${\cal{V}}_{W,N}^{p,q,r}$, $p,q,r\in \mathbb{Z}$, the family of neighborhoods
\begin{eqnarray*}
V(x)&=&\{y:(a_{-p}(\pi_1(x))\leq\pi_1(y)\leq a_1(\pi_1(x))\wedge a_{-q}(\pi_2(x))\leq\pi_2(y)\leq a_r(\pi_2(x)))\\
&& \ \ \ \ \vee (\pi_1(y)=\pi_1(x)\wedge a_1(\pi_2(y))\leq\pi_2(y)\leq a_r(\pi_2(y)))
\end{eqnarray*}
we obtain an analogous result to Lemma 3.1, if we assume condition $D''(u_n,{\cal{B}}_n,{\cal{V}}_{W,N}^{p,q,r})$.

In fact, we can decompose the event $\{\bigvee_{x\in B_n^{(s,t)}} Z(x)>u_n(\tau)\}$ in eight different ways corresponding to the different forms of starting from $x$ and getting around $\mathbb{Z}^2$, along the directions $\{(a_k(\pi_1(x)),0),k\in \{-1,1\}\}$ and $\{(0,a_k(\pi_2(x))),k\in \{-1,1\}\}$.
\end{rem}

\subsection*{Proof of Proposition \ref{proposition3.1}}

From Lemma \ref{lemma3.1} and the definition of the tail dependence coefficient $\lambda_n^{(\tau)}(V(x),x)$, it follows that
\begin{eqnarray*}
\lefteqn{P\left(\displaystyle\bigvee_{x\in B_n^{(s,t)}} Z(x)>u_n(\tau)\right)}\\
&=&\displaystyle\sum_{x\in B_n^{(s,t)}}P\left(Z(x)>u_n(\tau)\geq \displaystyle\bigvee_{y\in V(x)} Z(y)\right)+o\left(\frac{1}{k_n^2}\right)\\
&=&\overline{F}(u_n(\tau))\displaystyle\sum_{x\in B_n^{(s,t)}}P\left(\displaystyle\bigvee_{y\in V(x)} Z(y)\leq u_n(\tau) \mid Z(x)>u_n(\tau)\right)+o\left(\frac{1}{k_n^2}\right)\\
&=&{\overline {F}}(u_n)\displaystyle\sum_{x\in B_n^{(s,t)}}(1-\lambda_n^{(\tau)}(V(x),x))+o\left(\frac{1}{k_n^2}\right).
\end{eqnarray*}
Thus,
\begin{eqnarray*}
\lefteqn{\displaystyle\lim_{n\rightarrow\infty}P\left(\displaystyle\bigvee_{x\in A_n}Z(x)\leq u_n(\tau)\right)}\\
&=&\exp\left(-\displaystyle\lim_{n\rightarrow\infty}\displaystyle\sum_{B_n^{(s,t)}\in {\cal{B}}_n}P\left(\displaystyle\bigvee_{x\in B_n^{(s,t)}}Z(x)>u_n(\tau)\right)\right)\\
&=&\exp\left(-\displaystyle\lim_{n\rightarrow\infty}\displaystyle\sum_{B_n^{(s,t)}\in {\cal{B}}_n}\left(\overline{F}(u_n(\tau))\displaystyle\sum_{x\in B_n^{(s,t)}}(1-\lambda_n^{(\tau)}(V(x),x))+o\left(\frac{1}{k_n^2}\right)\right)\right)\\
&=&\exp\left(-\displaystyle\lim_{n\rightarrow\infty}\overline{F}(u_n(\tau))\displaystyle\sum_{B_n^{(s,t)}\in {\cal{B}}_n}\displaystyle\sum_{x\in B_n^{(s,t)}}(1-\lambda_n^{(\tau)}(V(x),x))\right)\\
&=&\exp\left(-\displaystyle\lim_{n\rightarrow\infty}\overline{F}(u_n(\tau))\displaystyle\sum_{x\in A_n}(1-\lambda_n^{(\tau)}(V(x),x))\right)\\
&=&\exp\left(-\tau\displaystyle\lim_{n\rightarrow\infty}\frac{1}{f(n)}\displaystyle\sum_{x\in A_n}(1-\lambda_n^{(\tau)}(V(x),x))\right),
\end{eqnarray*}
which concludes the proof.


\end{appendix}

\begin{appendix}

\section*{Appendix C: Proofs for Section 4}\setcounter{section}{0}

\refstepcounter{section}

\subsection*{Proof of Proposition \ref{proposition4.1}}

Since $u_n(\tau)=\frac{f(n)}{\tau}$ and $F(y)=\exp(-y^{-1})$, we have
\begin{eqnarray*}
\lefteqn{\displaystyle\lim_{n\rightarrow\infty}\lambda_n^{(\tau)}(V(x),x)}\\
&=&\displaystyle\lim_{n\rightarrow\infty}P\left(\displaystyle\bigvee_{y\in V(x)}Z(y)>\frac{f(n)}{\tau}\left|Z(y)>\frac{f(n)}{\tau}\right.\right)\\
&=&\displaystyle\lim_{n\rightarrow\infty}P\left(F\left(\displaystyle\bigvee_{y\in V(x)}Z(y)\right)>F\left(\frac{f(n)}{\tau}\right)\left| F(Z(x))>F\left(\frac{f(n)}{\tau}\right.\right)\right)\\
&=&\displaystyle\lim_{n\rightarrow\infty}P\left(\displaystyle\bigvee_{y\in V(x)}F(Z(y))>\exp\left(-\frac{\tau}{f(n)}\right)\left|F(Z(x))>\exp\left(-\frac{\tau}{f(n)}\right.\right)\right)\\
&=&\displaystyle\lim_{n\rightarrow\infty}P\left(\displaystyle\bigvee_{y\in V(x)}F(Z(y))>1-\frac{\tau}{f(n)}\left| F(Z(x))>1-\frac{\tau}{f(n)}\right.\right)\\
&=&\Lambda_U^{(V(x)\mid \{x\})}\\
&=&1+\epsilon_{V(x)}-\epsilon_{V(x)\bigcup\{x\}}
\end{eqnarray*}


\subsection*{Proof of Proposition \ref{proposition4.2}}

Results from the fact that $\{1-\lambda_n^{(\tau)}(V(x),x)\}_{n\geq 1}$ converges to $\theta_{\{x\}}$, uniformly in $x$, since this implies that for every $\epsilon>0$, there exists a natural number $m=m(\epsilon)$ such that, for every $n>m$,
\begin{eqnarray*}
\frac{1}{f(n)}\displaystyle\sum_{x\in A_n}(1-\lambda_n^{(\tau)}(V(x),x))&\leq&\frac{1}{f(n)}\displaystyle\sum_{x\in A_n}(\theta_{\{x\}}+\epsilon)\\
&=&\left(\frac{1}{f(n)}\sum_{x\in A_n}\theta_{\{x\}}\right)+\epsilon
\end{eqnarray*}
and
\begin{eqnarray*}
\frac{1}{f(n)}\displaystyle\sum_{x\in A_n}(1-\lambda_n^{(\tau)}(V(x),x))&\geq&\frac{1}{f(n)}\displaystyle\sum_{x\in A_n}(\theta_{\{x\}}-\epsilon)\\
&=&\left(\frac{1}{f(n)}\sum_{x\in A_n}\theta_{\{x\}}\right)-\epsilon.
\end{eqnarray*}


\end{appendix}

\end{document}